\documentclass{article}
\usepackage{amssymb,amsmath, tikz}
\usepackage{graphics}
\input xy
\xyoption{all}
\usepackage{accents}
\usepackage{braket}
\usepackage{hyperref}
\newlength{\dhatheight}

\newcommand\phantomarrow[2]{%
  \setbox0=\hbox{$\displaystyle #1\to$}%
  \hbox to \wd0{%
    $#2\mapstochar
     \cleaders\hbox{$\mkern-1mu\relbar\mkern-3mu$}\hfill
     \mkern-7mu\rightarrow$}%
  \,}

\begin{document}

\title{Induced $C^*$-complexes in metaplectic geometry}

\author{Svatopluk Kr\'ysl\footnote{{\it E-mail address}: 
Svatopluk.Krysl@mff.cuni.cz}\\ {\it \small Charles University, Faculty of 
Mathematics and  Physics,  Czechia}
\thanks{The author thanks for financial supports from the founding 
No. 17-01171S granted by the Czech Science Foundation. We thank to  the 
anonymous reviewer for his  
 comments and suggestions.}}

\maketitle \noindent
\centerline{\large\bf Abstract}  

For a symplectic manifold admitting a metaplectic structure and for a Kuiper map, we 
construct a complex of differential operators acting on exterior differential 
forms with values in the   dual of  Kostant's symplectic spinor bundle. 
Defining a Hilbert $C^*$-structure on this  bundle for a suitable $C^*$-algebra,  
we  obtain an elliptic $C^*$-complex in the sense of Mishchenko--Fomenko. Its cohomology groups 
appear to be finitely generated projective Hilbert $C^*$-modules. The paper can 
serve as a guide for handling of  differential complexes and  
PDEs on Hilbert bundles.

\section{Introduction}
   
An important assumption is made in 
classical works on the  Hodge theory of elliptic complexes beside the compactness of the base manifold.
Namely, the bundles in the complexes are supposed to be of finite rank. (See Hodge 
\cite{Hodge} and Palais \cite{Palais}.) Under these two assumptions, 
elements of the complexes can be decomposed into closed, co-closed 
and harmonic parts uniquely. It follows that cohomology groups of the 
complexes are  bijectively representable 
by  harmonic elements. Moreover, they are finite 
dimensional vector spaces.
The study of elliptic complexes has a long history which goes back to 
W. Hodge  and G. de Rham. See  Maurin \cite{Maurin} and also Klein 
\cite{Klein} for references and for the (pre-)history of this subject.
Elliptic complexes  are  studied also in areas `extraterritorial' to 
differential geometry and global analysis or in cases 
violating the compactness. See Tsai at 
al. \cite{Tsai}, Schmid, Vilonen \cite{SchmidVilonen}, Li \cite{Li},  
Nekov\'a\v{r}, Scholl 
\cite{Nekovar},  Hain \cite{Hain} and Albin et al. \cite{ALMP} for  recent 
contributions.

However, elliptic complexes can be studied in the case of 
infinite rank bundles as well. 
Troitsky in \cite{Troitsky} and  Schick, subsequently, in \cite{Schick} elaborate a theory for 
elliptic complexes on   finitely generated 
projective Hilbert $C^*$-bundles, which is based on the index theory of 
Mishchenko and Fomenko described in \cite{Mishchenko} and \cite{FM}.
In \cite{Troitsky} and \cite{Schick}, $C^*$-compact perturbations of the
complexes' differentials are allowed. For the infinite rank bundles, one cannot 
expect that, in general,
the space of harmonic elements represents the cohomology groups
homeomorphically (quotient topology) and linearly (quotient projections). 
Therefore, it makes sense to find conditions when this happens. This was 
investigated by the author in the  past years (see \cite{KryslAGAG2}, \cite{KryslJGP2} or 
an overview in \cite{KryslHabil}).

 In the connection to  sheaf cohomology, Banach and Fr\'echet bundles are 
studied in the  papers of Illusie 
\cite{Illusie} and R\"ohrl \cite{Rohrl}. To present a sample of the broad 
context 
(foremost connected to $C^*$-algebras, $K$-theory and homological algebra) in which 
such bundles are considered, let us mention the papers of Maeda, Rosenberg 
\cite{MaedaRosenberg}, Freed, Lott \cite{Freed}, Larra\'in-Hubach 
\cite{Larrain}, the author \cite{KryslJGP2} and  Fathizadeh, Gabriel 
\cite{FathiGabr}. There are also  works  
in which holomorphic Banach bundles and their sheaf cohomology groups are 
treated. See  Lempert \cite{Lempert} and Kim \cite{Kim}. A further reason  for 
an investigation of these complexes might originate in BRST-theories 
(\cite{HenTeitel}), which we do not touch here explicitly, although we give a 
modest interpretation of the introduced  structures in the realm of quantum 
theory.

The aim of our paper is to study elliptic complexes for the case of infinite 
rank Hilbert $C^*$-bundles which are 
induced by Lie group representations, and to show as well 
how the theory of connections and  appropriately generalized 
elliptic complexes  apply in this situation. In order to make the considered
situation not too general, 
we choose a specific Lie group representation (the Segal--Shale--Weil representation of the metaplectic group). 
We hope that the reader may follow the text easier. Our further aim  
is to find  examples for the theory of Mishchenko and Fomenko 
(\cite{FM}) that would not be trivial or covered by the theory developed around 
the classical Atiyah--Singer index theorem, and to show how the Hodge theory for certain $C^*$-bundles
 (derived  in \cite{KryslAGAG2} and \cite{KryslJGP2}) make us able to describe 
the cohomology of elliptic complexes 
on Hilbert bundles easily.

Let $(V,\omega_0)$ be a symplectic vector space and $L$ be a Lagrangian 
subspace of it. The Segal--Shale--Weil representation is a non-trivial representation of 
the connected double cover $\widetilde{G}$ of the symplectic group  on a 
complex Hilbert space by unitary operators. The Hilbert space is formed by the 
$L^2$-functions defined on $L.$
For references to this representation see, e.g., Weil 
\cite{Weil}, Borel, Wallach \cite{BorelWallach}, Wallach \cite{Wallach}  and 
Folland \cite{Folland}. Let us denote the dual of the Segal--Shale--Weil 
representation  by $\sigma$ and the continuous dual of the Hilbert space by $H.$ We call the resulting representation   
the {\it oscillator representation}.
The double cover, the  metaplectic group, acts by
$\sigma_k: \widetilde{G} \to \mbox{Aut}(H^k)$ on the tensor products 
$H^k=\bigwedge^k V^*\otimes H$ in a natural way.
Let $CH$ denote the algebra of compact operators on 
the vector space carrier of the Segal--Shale--Weil representation, i.e., on $L^2(L).$  
On the vector space space $CH,$ we consider the representation, denoted by 
$\rho,$ that is given by the conjugation by the  
Segal--Shale--Weil representation.

We introduce a right $CH$-module structure $H \times CH \to H$ on $H$
and a  $CH$-product $H\times H \to CH$ obtaining a Hilbert $CH$-module -- a 
mixture of a module and of a normed vector space as
described by Paschke   \cite{Paschke} and  Rieffel  \cite{Rieffel}. 
The 
spaces $H^k,$ $k\geq 1,$ 
are made Hilbert $CH$-modules as well. 
Inducing $(\sigma_k,H^k)$ to a principal $\widetilde{G}$-bundle over a
symplectic manifold $M$ and in particular, to the so-called metaplectic 
bundle of a metaplectic structure (Kostant 
\cite{Kostant}), one 
gets Banach bundles $\mathcal{H}^k$ (Sect. 2).
Further,  we associate $(\rho, CH)$ 
to the metaplectic bundle, getting a 
continuous Banach bundle. We  denote it by $\mathcal{CH}$ and call it the bundle of measurements.
Then we define certain intermediate maps $\mathcal{H} \times_{M} 
\mathcal{CH} \to \mathcal{H}$ and $\mathcal{H} \times_{M} \mathcal{H} \to 
\mathcal{CH},$ called pointwise structures, which are families of 
$C^*$-algebra actions and of $C^*$-algebra valued products.
We are not aware that any substantial modification was 
done regarding the assumptions in the theory of Mishchenko--Fomenko from \cite{FM}.
Thus, we adapt the  pointwise structures to get the honest required ones.  
Recall that we have to  define maps $\mathcal{H} \times CH \to \mathcal{H}$ and 
$\mathcal{H} \times_M \mathcal{H} \to CH$ instead of $\mathcal{H} \times_M 
\mathcal{CH} \to \mathcal{H}$  and
$\mathcal{H}\times_M \mathcal{H} \to \mathcal{CH},$ respectively.
This can be assured by choosing a trivialization $J$ of $\mathcal{H},$ called a 
Kuiper map, which is possible in the case of a unitary representation.
Nevertheless, let us note
that we can  choose the  maps 
 $\mathcal{H} \times CH \to \mathcal{H}$ and $\mathcal{H} \times_M \mathcal{H} 
\to CH$ by using a Kuiper map without defining the 
pointwise structures. We do not do so  and 
work with the pointwise maps as far as possible instead. Both the $CH$-action 
and the 
$CH$-valued fiber product on $\mathcal{H}^{\bullet}$ are defined in Sect. 4.2.

We focus our attention on the continuity and differentiability of 
bundles induced by unitary or isometric representations of Lie groups  
on infinite dimensional Hilbert or Banach spaces, respectively.
We show that one may  speak  about $C^{\infty}$-differentiable structures  
when a particular smooth principal bundle atlas is chosen, but only on $C^0$-differentiable structures for general representations.
 We treat connections, vertical lifts, connectors and induced covariant 
derivatives in the context of Banach bundles, showing   that the general
  approach to connections as given,  e.g., in Kol\'a\v{r}, Michor, 
Slov\'ak \cite{KMS} is possible to 
apply in the case of  infinite bundles as well.
Since the oscillator representation 
is known  to be  not differentiable, we do not get 
a $C^{\infty}$-structure that is induced canonically. As already mentioned, it is possible to 
choose a Kuiper map. The Kuiper structure on $\mathcal{H}$ is defined, in a non-canonical way,  as the maximal
smooth atlas containing the Kuiper map.
Some facts from Lie groups' representations based on Bruhat \cite{Bruhat} ($l=\infty$) and 
Goodman \cite{Goodman} ($l \in \mathbb{N}_0$) provide 
us with further $C^l$-structures as well.   We define them in 
Sect. 3.2.  However, only the Kuiper structure seems to have 
some purpose in the theory of $C^*$-indices and consequently, for
$C^*$-elliptic complexes.

Let us mention that we don't use other algebras than the $C^*$-algebras  
for several 
reasons. The first one is that we would like to stay close to the approach of 
quantum physics, for which  $C^*$-algebras are fundamental objects until now. 
The second one is more pragmatic. As far as we know, the Hodge 
theory has not yet been developed for bundles over other algebras in a 
satisfactory form. Let us mention that a different approach might be represented 
possibly
by considering bundles over Fr\'echet $*$-algebras (see Keyl et al. 
\cite{Keyl} for these algebras and their connection to quantum physics). This 
approach eliminates problems with the smoothness of the representations.

Further, we give a straightforward proof of the triviality of 
the bundle of measurements (Sect. 3.3), which is not easy to find in the
literature, but which is colloquial if one is acquainted with certain constructions in
non-abelian sheaf cohomology adapted for Banach bundles. 

  In the last section of the paper, we use the Kuiper map of the oscillator 
bundle $\mathcal{H}$ to define a 
product type  connection $\Phi$  
which is flat by construction and induces a cochain 
complex  $d_{\bullet}^{\Phi}=(\Gamma(\mathcal{H}^k), d_k^{\Phi})_{k\in 
\mathbb{Z}}$ of exterior covariant derivatives acting between sections of 
$\mathcal{H}^k.$
When the symplectic manifold is compact, we prove the central result of the 
paper (Theorem 18) that the cohomology groups of $d_{\bullet}^{\Phi}$ are 
isomorphic to 
the tensor product of the de Rham groups and $H.$ 
 We  show also how to use the Hodge 
theory 
for  projective and finitely generated Hilbert $CH$-bundles in the case of  
$d_{\bullet}^{\Phi},$ that is elliptic, for to produce a substantially shorter proof of Thm. 18.

\bigskip

In order to specify our terminology, we start with a preamble to 
which the reader may come back during the reading.

\bigskip

{\bf \large Preamble}

\begin{itemize}  
\item[a1)] Let $X$ be a topological space. The symbol $\mbox{Aut}(X)$ denotes 
the set of all homeomorphisms of $X.$ 
\item[a2)] All vector spaces are considered to be  Fr\'echet. Hilbert spaces 
are considered with the norm topologies induced by  
scalar products and the scalar products are anti-linear in the left variable 
(physicists' convention).
\item[a3)] If $X$ is a vector space, maps in $\mbox{Aut}(X)$ are additionally 
supposed to be linear.
By $\mbox{End}(X),$ we understand continuous linear maps $X \to X.$ In 
particular, $\mbox{End}(X)=B(X)$ where $B(X)$ are bounded operators on $X.$
By  a functional on a Banach space, we mean a continuous linear 
map of this space into the field of complex numbers
\item[b)] By a $C^l$-differentiable structure on a manifold, we understand a
maximal  $C^l$-atlas with respect to the inclusion. A 
$C^{\infty}$-differentiable structure is called a smooth structure. 
\item[c)] All manifolds are $C^0$-differentiable (topological) Banach 
manifolds and all base spaces 
of fiber bundles are finite dimensional smooth manifolds. All 
bundles are locally trivial.
\item[d)] Representations of Lie groups are considered 
continuous in the classical representation theory sense. For a homomorphism 
$\sigma: G \to 
\mbox{Aut}(X),$ of a Lie group $G$ into the linear homeomorphisms of a normed 
vector space $X$ to be a representation  means 
 that  $\widetilde{\sigma}: G \times X \to X$ given by $\widetilde{\sigma}: G 
\times 
H \ni (g,f) \mapsto \sigma(g)f \in X$  is continuous.\footnote{If more details 
needed, see, e.g., Knapp 
\cite{KnappOverview}.} 
\item[e)] When an index of a labeled  object or a labeled morphism exceeds its 
allowed range, we consider it as zero. 
\end{itemize}

\section{$\widetilde{G}$-module and $CH$-module structures}

In this section, we introduce the higher oscillator representation of the metaplectic 
group and the Hilbert $CH$-module structure on the underlying space of this representation,
 and we prove some properties of these structures.

Let $(V,\omega_0)$ be a real symplectic vector space of dimension 
$2n \geq 1$ and  $G$ denote the  symplectic group $Sp(V,\omega_0).$ 
Let us choose  a compatible positive complex structure $J_0$ on $V.$ By a 
compatible positive complex structure,
we mean an $\mathbb{R}$-linear map $J_0:V \to V$ such that $J_0^2 = -\mbox{Id}_V$  and 
such that $g_0$ defined by 
$g_0(u,v) = \omega_0(J_0u,v),$ $u,v \in V,$ is a positive definite bilinear 
form.
Denoting the orthogonal group of $(V,g_0)$ by $O(V,g_0),$  we set $K = 
O(V,g_0)\cap Sp(V,\omega_0)$ for the 
unitary group of the hermitian vector space $(V,J_0,g_0).$

The symplectic group  is known to retract smoothly onto $K$ 
(see, e.g., \cite{Folland}) the first fundamental group of which is 
isomorphic to $\mathbb{Z}.$
By the theory of covering spaces, there exists a unique smooth connected 
$2$-fold covering  
$\lambda:\tilde{G} \to G$ of $G$ up to a smooth covering isomorphism.
 Choosing an element in the $\lambda$-preimage of the neutral element of $G$ 
and declaring it as the neutral element of $\widetilde{G},$
define a unique Lie group
structure  on $\widetilde{G}.$ 
This is the so called {\it metaplectic group}, denoted 
by $Mp(V,\omega_0).$  Further, we set $\widetilde{K} = \lambda^{-1}(K)$ for 
the corresponding double cover of the unitary group.

\bigskip

{\bf Basics on oscillator representation:} After Shale \cite{Shale} 
and Weil \cite{Weil}, it is known that there exists a specific 
non-trivial unitary 
representation
of $\widetilde{G},$  called metaplectic, Shale--Weil, Segal--Shale--Weil, 
oscillator or symplectic spinor representation. It is called a 
spinor representation probably because it
does not descend to a true representation of the symplectic group.   
Since 
$\widetilde{G}$ is non-compact,  
the space carrying this representation is an infinite dimensional Hilbert 
space.
The representation intertwines so-called twisted Schr\"{o}dinger representations 
and 
can be made a representation (of the metaplectic group) on the Hilbert space 
$L^2(L)$ of  
complex valued square Lebesgue integrable functions defined on a Lagrangian 
subspace
$L \subseteq (V,\omega_0).$ The Lebesgue measure is determined by the norm on 
$L$ induced
by the inner product $g_0$ restricted to $L \times L.$ For the uniqueness of the 
representation,
we refer to Wallach \cite{Wallach}, p. 224. We denote the Hilbert 
space 
$L^2(L)$ by $\check{H}$ and  fix a Segal--Shale--Weil representation, denoting 
it by
$\check{\sigma}: \widetilde{G} \to U(L^2(L))$ where $U(X)$ 
denotes the space of unitary operators on a Hilbert space $X.$

\bigskip

{\bf Notation for dualities $\&$ musical isomorphisms:} Let $(\check{H}, 
(,)_{\check{H}})$ be a  complex Hilbert space, and let us denote its
continuous dual by $H.$
For each $f \in  H,$ we denote its dual  by $f^{\sharp},$ which is the unique 
element in  $\check{H}$ such that
$f(v)=(f^{\sharp},v)_{\check{H}}$ for any $v \in \check{H}.$ By the Riesz's 
representation theorem for Hilbert spaces, the element $f^{\sharp}$
exists also for non-separable Hilbert spaces. 
For any  $v \in \check{H},$ we denote by $v^{\flat}$ the element of $H$ defined 
by
$v^{\flat}(w)=(v,w)_{\check{H}},$ $v,w \in \check{H}.$ We set $(f,g)_H = 
(g^{\sharp},f^{\sharp})_{\check{H}}$ in order the product be anti-linear in the 
left variable.
With this notation, $(H,(,)_H)$ is a complex Hilbert space and $\sharp: 
\check{H} \to H$ and $\flat:H \to \check{H}$ are 
anti-unitary maps, i.e., they are complex anti-linear and   intertwine the 
corresponding scalar products, i.e.,
$(v^{\flat}, w^{\flat})_H = (w,v)_{\check{H}},$ $(f^{\sharp}, \, 
g^{\sharp})_{\check{H}} = (g,f)_H$.
For any pair $u \in \check{H}$ and $v \in H,$ we denote by $u \otimes v$ the 
endomorphism $u \otimes v: w \ni \check{H} \mapsto v(w)u \in \check{H}.$
For any linear continuous map $A:X \to Y$ between Hilbert spaces $X$ and $Y,$ 
$A^*:Y \to X$ denotes its adjoint.

\bigskip

Because we would like to consider right $C^*$-modules, it is more convenient to 
use the dual of the Segal--Shale--Weil representation. We denote it 
by $\sigma$ and call it the {\it oscillator representation}. Using the  
notation defined, we  have a homomorphism 
$$\sigma: \widetilde{G}  \to U(H)$$ The homomorphism is a 
representation, i.e., it is continuous in the sense described in 
Preamble, Par. d.

Let us notice that the oscillator representation decomposes into two 
irreducible subrepresentations, the space of functionals on  even and that 
of odd functions in $\check{H}.$  
See   Wallach \cite{Wallach} and Folland \cite{Folland}
for more information on $Mp(V,\omega_0)$ and the Segal--Shale--Weil 
representation.

Since $V$ is equipped with $g_0$ and $\omega_0$ (the latter gives an 
orientation 
to $V$ which is induced by $\omega_0^{\wedge n}$),
 we get a canonical scalar product (Hodge type product) on $\bigwedge^k V^*,$ 
$k=0,\ldots, 2n.$ We denote it by $g_0$ as well.
For $k=0,\ldots, 2n,$ the tensor product
$H^k = \bigwedge^k V^* \otimes H$  is equipped with the norm topology induced by 
the canonical Hilbert space inner 
product on the tensor product of Hilbert spaces.  Let 
${\lambda^*}^{\wedge k}$ denote the $k$-th wedge product of
the dual of  the representation (and covering) $\lambda:\widetilde{G}\to 
Sp(V,\omega_0) \subseteq \textrm{Aut}(V).$  
 The tensor product 
representations 
$\sigma^k:\widetilde{G} \to \mbox{Aut}(H^k)$ are given by  
$$\sigma^k(g)(\alpha \otimes f) = {\lambda^*}^{\wedge k}(g)\alpha\otimes 
\sigma(g)f$$ where
$g\in \widetilde{G},$ $\alpha \in \bigwedge^k V^*$ and  $f\in H.$ 
The above prescription for $\sigma^k$ is extended linearly to 
non-homogeneous 
elements.
Let us remark that $H^0 = H$ is the representation space of the oscillator 
representation $\sigma =\sigma^0.$ 
We set 
$H^{\bullet} = \bigoplus_{k=0}^{2n}H^k 
=\bigoplus_{k=0}^{2n} \bigwedge^k V^* \otimes H.$ 
In a parallel to ordinary spinors, we call 
$H^{\bullet}=\bigoplus_{k=0}^{2n} H^k$ the {\it higher oscillator module} and 
the representation
$\sigma^{\bullet} = \bigoplus_{k=0}^{2n}\sigma^k$ the {\it higher oscillator 
representation}. 

\bigskip

Let us denote the space of compact operators on $\check{H}$ equipped with the 
operator norm by $CH,$ and consider the representation $$\rho: 
\widetilde{G} \to 
\mbox{Aut}(CH)$$ of the metaplectic group $\widetilde{G}$ on $CH$
 given by $\rho(g)a = \check{\sigma}(g) \circ a \circ \check{\sigma}(g)^{-1},$ 
where $g \in 
\widetilde{G}$ and $a \in CH.$ 
The representation is well established since $CH$ is a closed two-sided
ideal in the space $\textrm{End}(\check{H})$ of bounded operators on 
$\check{H}.$
 
\bigskip

With the adjoint of a Hilbert space map as the $*$-operation, 
$(CH,||\,||_{CH},*)$ becomes a $C^*$-algebra, where $||\,||_{CH}$ is the 
restriction 
of the operator norm on $B(\check{H})$ to $CH.$

\bigskip

{\bf Definition:}
The {\it right  action} $H^{\bullet} \times CH \to H^{\bullet}$
of $CH$ on $H^{\bullet}$ is  defined on homogeneous elements by  $$(\alpha 
\otimes u) \cdot a = \alpha 
\otimes (u \circ a)$$ where $\alpha \otimes u \in 
\bigwedge^{\bullet} V^* \otimes H$ and $a\in CH.$
It is extended  linearly to all elements of $H^{\bullet}.$  
The {\it $CH$-product} $(,): H^{\bullet}\times H^{\bullet} \to CH$ is given  by 
$$(\alpha \otimes u, \beta \otimes v) = 
g_0(\alpha,\beta)  u^{\sharp} \otimes v$$ where
$\alpha, \beta \in \bigwedge^{\bullet} V^*$ and $u,v \in H.$  
The product is extended by linearity. It maps $H^{\bullet} \times H^{\bullet}$ 
into finite rank operators, especially, into $CH.$
We denote the induced ($C^*$-module) norm on $H^{\bullet}$ by $||\,||.$ It is 
given by $||f|| = 
\sqrt{||(f, f)||_{CH}},$ $f\in H^{\bullet}.$

\bigskip

{\bf Lemma 1:} For $k=0,\ldots, 2n$ and $\alpha \otimes f \in H^k$ 
$$||\alpha \otimes f|| = \sqrt{g_0(\alpha,\alpha)}||f||_H$$ In 
particular, the norm $||\,||$ coincides on homogeneous elements with the norm induced by the inner 
product on $H^k.$

{\it Proof.} The appearance of the first factor, $\sqrt{g_0(\alpha, \alpha)},$ 
is easy to show. For the second one, we compute
the operator norm $||\,||_{CH}$ of the map $f^{\sharp} \otimes f: \check{H} \ni 
k 
\mapsto f(k)f^{\sharp},$ where $f \in H.$ On one hand,
 we have $||(f,f)(k)||_{\check{H}}^2=||f^{\sharp}f(k)||^2_{\check{H}}$ $=|f(k)|^2 
(f^{\sharp},f^{\sharp})_{\check{H}}=
|f(k)|^2(f,f)_H \leq $ $||f||^2_H$
$||k||_{\check{H}}^2$ 
$||f||_H^2.$ 
Consequently, 
$||(f,f)||_{CH}  \leq ||f||_H^2.$
On the other hand,
for $k=f^{\sharp}\neq 0,$ we get  
$||f^{\sharp}f(f^{\sharp})||_{\check{H}}^2/$ $/||f^{\sharp}||_{\check{H}}^2
= ||f||^2_H(f,f)_H/||f||^2_H=||f||^2_H.$ Thus, 
$||(f,f)||_{CH}=||f||_H^2=(f,f)_H.$
For the norm  induced by the $CH$-product on  $H$, we obtain 
$||f||=\sqrt{||(f,f)||_{CH}}=||f||_H,$
and consequently, $||f||=\sqrt{(f,f)_H}.$
\hfill\(\Box\)

\bigskip

{\bf Lemma 2:} The  pair $(H^{\bullet}, (,))$ is a finitely generated 
projective Hilbert $CH$-module.

{\it Proof.}  
It is obvious that $H^{\bullet}$ is a right $CH$-module.
Proving that $(,)$ is hermitian symmetric and $CH$-sesquilinear, is 
a routine 
check. Indeed, for $a\in CH,$ $\alpha \otimes f, \beta \otimes h \in H^k,$ and 
$u\in 
\check{H}$ 
 \begin{align*}
  ((\alpha \otimes f) \cdot a, \beta \otimes h)(u)&= g_0(\alpha, 
\beta)\left((f\cdot a)^{\sharp} \otimes h \right)(u)
=g_0(\alpha, \beta) \left((f \circ a)^{\sharp} \otimes h \right)(u)\\
&=g_0(\alpha,\beta) \left(a^*(f^{\sharp}) \otimes h \right)(u)=g_0(\alpha, 
\beta)  h(u) a^*(f^{\sharp})\\
&=a^*\left(g_0(\alpha, \beta)h(u) f^{\sharp}\right)
=a^*\left(g_0(\alpha,\beta) \left(f^{\sharp} \otimes h\right) (u) \right)\\
&=\left(a^*\circ \left(g_0(\alpha, \beta) \left(f^{\sharp} \otimes h 
\right)\right)\right) (u)\\ 
&= \left(a^* \circ (\alpha \otimes f, \beta \otimes 
h)\right) (u)
\end{align*}
Thus, $(,)$ is anti-linear over $CH$ in the left variable. The hermitian 
symmetry 
is proved similarly. The right $CH$-linearity follows from the hermitian 
symmetry and
left $CH$-linearity.
 
The projectivity follows from the Magajna theorem (Theorem 1, 
Magajna \cite{Magajna}).

For the existence of a finite set of generators, let us choose a basis 
$\mathfrak{C}$ of $\bigwedge^{\bullet}V^*$ and take a non-zero vector
$w \in H.$ For any $v \in H^{\bullet},$ there exist scalars 
$\lambda_{\mathfrak{a}} \in \mathbb{C}$ and vectors $v_{\mathfrak{a}} \in H,$ 
$\mathfrak{a} \in \mathfrak{C},$ such that
$v= \sum_{\mathfrak{a} \in \mathfrak{C}} \lambda_{\mathfrak{a}} \mathfrak{a} 
\otimes v_{\mathfrak{a}}.$   
For each $v' \in H,$ let us set
$p_{v',w}(k)=  \frac{v'(k)}{||w||_H^2} w^{\sharp},$  $k \in \check{H}.$
Then $v=\sum_{\mathfrak{a} \in \mathfrak{C}} (\mathfrak{a} \otimes w) \cdot 
(\lambda_{\mathfrak{a}} p_{v_{\mathfrak{a}},w}).$
Thus, $\{\mathfrak{a} \otimes w|\, \mathfrak{a} \in \mathfrak{C}\}$ represents a 
finite set of generators of 
$H^{\bullet}$.

Since the $C^*$-module norm $||\,||$ restricted to $H^0$ is the dual of the 
Hilbert
norm on  $\check{H}$ (Lemma 1), we see that $H^{\bullet}$ is complete.
\hfill\(\Box\)

\bigskip

{\bf Remark:} A similar proof  can be found in  \cite{KryslDGA} 
where, though, a determination of the inducing scalar product for the  
$C^*$-norm is 
omitted. The above construction of a Hilbert $CH$-module 
structure can be done for an arbitrary Hilbert space.

 Note that in the bra--ket convention, $p_{v',w} = 
\frac{\ket{w}\bra{v'}}{\braket{w|w}} = \frac{\ket{w}\bra{v'}}{||w||^2}.$ 

\bigskip

{\bf Definition:} For $a,b \in CH,$ $u \in \check{H}$ and $h \in H,$ let us 
set

\begin{tabular}{lll}
                $\circ: CH \times CH \to CH$ & $\circ 
(a,b)= a \circ b$ & (composition)\\
                $\otimes: \check{H} \times H \to CH$ & 
$ \otimes(u,h) = u \otimes h$ & (tensor product) \\
                $ev: H \times \check{H} \to 
\mathbb{C}$ & $ev(h,u) = h(u)$ & (evaluation) \\
\end{tabular}

\bigskip

{\bf Lemma 3:} The $CH$-action, the composition, the tensor 
product and the evaluation 
are $\widetilde{G}$-equivariant. The $CH$-product  $(,): H \times H \to CH$ and 
the induced norm $||\,||:H^{\bullet} \to \mathbb{R}$ are
$\widetilde{K}$-equivariant.
For any $g \in \widetilde{G},$ $v \in \check{H}$ and $h \in H$  
\begin{eqnarray} \label{formule}
\rho(g) (v \otimes h) = \check{\sigma}(g)v \otimes \sigma(g)h
\end{eqnarray}

{\it Proof.} For $g \in \widetilde{G},$ $v, w \in \check{H}$ and $h \in H$ 
\begin{align*}
\rho(g)(v \otimes h)(w) &= \left(\check{\sigma}(g)\circ (v \otimes h)\circ 
\check{\sigma}(g)^{-1}\right)(w)\\
&=\check{\sigma}(g)(v)h(\check{\sigma}(g)^{-1}w) \\
&= (\check{\sigma}(g)v \otimes \sigma(g)h)(w)
\end{align*}
obtaining formula   (\ref{formule}).

  (\ref{formule}) is used  to prove the $\widetilde{G}$-equivariance 
 of the tensor product readily.
We check the  $\widetilde{G}$-equivariance of the action $H^{\bullet}\times CH 
\to H^{\bullet}$ and the $\widetilde{K}$-equivariance of the 
$CH$-product. 
For $g \in \widetilde{G},$ $f =\alpha \otimes u \in H^k,$ $k=0,\ldots, 2n,$ and 
$a \in CH$  
\begin{align*}
(\sigma^k(g)f) \cdot (\rho(g)a) &= \left(\lambda^{* \wedge k}(g)\alpha \otimes 
\sigma(g) u\right) \cdot (\check{\sigma}(g) \circ a \circ 
\check{\sigma}^{-1}(g))\\
& = \lambda^{* \wedge k}(g)\alpha \otimes \left(\sigma(g) u \circ 
(\check{\sigma}(g) \circ a \circ \check{\sigma}^{-1}(g))\right)\\
&  = \lambda^{* \wedge k}(g)\alpha \otimes \left((u  \circ 
\check{\sigma}^{-1}(g))\circ (\check{\sigma}(g) \circ a \circ 
\check{\sigma}^{-1}(g))\right)\\
&  =\lambda^{* \wedge k}(g)\alpha \otimes (u \circ a \circ 
\check{\sigma}^{-1}(g))\\
&  =\lambda^{* \wedge k}(g)\alpha \otimes \sigma(g)(u  \circ a) = 
\sigma^k(g)(\alpha \otimes u \cdot a) = \sigma^k(g)(f \cdot a)
\end{align*}
The equivariance on non-homogeneous elements follows from linearity.

Let $g \in \widetilde{K}$ and $h = \beta \otimes v \in H^k.$ Since $g \in 
\widetilde{K},$ we  have \\
$g_0(\lambda^{* \wedge k}(g)\alpha, \lambda^{* \wedge k}(g)\beta)= g_0(\alpha, 
\beta).$ Further
\begin{align*}
(\sigma^k(g)f, \sigma^k(g)h) &= g_0(\lambda^{* \wedge k}(g)\alpha,\lambda^{* 
\wedge k}(g)\beta)\left((\sigma(g)u)^{\sharp} \otimes (\sigma(g)v)\right)\\
& = g_0(\alpha, \beta) \check{\sigma}(g)u^{\sharp} \otimes \sigma(g)v \\
& = g_0(\alpha, \beta) \rho(g)(u^{\sharp} \otimes v) = \rho(g) g_0(\alpha, 
\beta) (u^{\sharp} \otimes v) \\
& = \rho(g)(\alpha \otimes u, \beta \otimes v) = \rho(g)(f,h)
\end{align*} 
where formula (\ref{formule}) is used in the second row.
For the norm, consider $||\sigma^k(g)f||^2=||(\sigma^k(g)f,\sigma^k(g)f)||_{CH} 
=||\rho(g)(f,f)||_{CH}=||(f,f)||_{CH}=||f||^2,$ $g\in \widetilde{K},$ 
since $\rho(g)$ is  isometric.

$\widetilde{G}$-equivariances of  $\circ$ and $ev$ can be proved in a similar 
way.
\hfill\(\Box\)

\section{Geometric lifts of $\widetilde{G}$- and $CH$-module structures}

In this chapter, we introduce the higher oscillator bundle and the bundle of 
measurements. We investigate continuity and smoothness of canonical 
atlases of these bundles. For this, we examine the
Segal--Shale--Weil representation from the analytic point of view.

\bigskip

For a symplectic manifold  $(M^{2n},\omega),$ we consider  
$$\mathcal{Q} = \{ A: V \to T_m^*M| \, \omega_m(Au, Av) = \omega_0(u,v) \mbox{ for any } u,v \in V, m \in M \}$$
together with the bundle projection $\pi_Q(A)=m$ if and only if  $A:V \to T_m^*M.$
For the topology on $\mathcal{Q},$ we take the final one with respect to  
inverses of all bundle maps for the projection  $\pi_Q.$  
The right action of $G=Sp(V,\omega_0)$ on $\mathcal{Q}$ is given by the composition from the right, i.e., 
$(A, g) \mapsto A \circ g$ for $A\in \mathcal{Q}$ and $g \in G.$ This make $\pi_Q: \mathcal{Q} \to M$ a principal
$G$-bundle.

\bigskip

Let us consider
pairs $(\pi_P,\Lambda),$ where $\pi_P:\mathcal{P}\to M$ 
is a principal $\widetilde{G}$-bundle  over $M$  
and $\Lambda: \mathcal{P} \to \mathcal{Q}$
is a smooth map. Such a pair is 
called a {\it metaplectic structure} 
if $\Lambda$ intertwines the actions of $G$ and $\widetilde{G}$ on 
$\mathcal{Q}$ and  $\mathcal{P},$ respectively,  with respect to the homomorphism
$\lambda: \widetilde{G} \to G$.  
In other words, the  following diagram has to  commute.
$$\begin{xy}\xymatrix{
\mathcal{P} \times \widetilde{G} \ar[dd]^{\Lambda\times \lambda} \ar[r]&   
\mathcal{P} \ar[dd]^{\Lambda} \ar[dr]^{\pi_P} &\\
                                                            & &M\\
\mathcal{Q} \times G \ar[r]   & \mathcal{Q} \ar[ur]_{\pi_Q} }\end{xy}$$
(The horizontal arrows  represent the right actions of the groups 
$\widetilde{G}$ and 
$G$ on the   total spaces $\mathcal{P}$ and $\mathcal{Q},$ 
respectively.) The bundle $\mathcal{P}$ is called a {\it metaplectic bundle}.
 
 The bundle $\mathcal{Q}$ is  isomorphic to the quotient of $\mathcal{P}$ by the 
action of the $2$-element group $\mbox{Ker} \, \lambda.$ 
Further, we may consider a  version of a metaplectic structure over the complex numbers 
$\mathbb{C}$ as well,
getting a so-called  $Mp^c$-structure. Notice  that a metaplectic structure is formally similar to that of a spin structure.

\bigskip

{\bf Remark:}
 Let us mention that a  metaplectic 
structure exists on a 
symplectic manifold $(M,\omega)$ if and only if 
the second Stiefel--Whitney class of $TM$ vanishes. Their set modulo an 
appropriate equivalence relation is isomorphic to the singular
cohomology group $H^1(M,\mathbb{Z})$ as a set.  See Kostant  
\cite{Kostant}.
Note that the $Mp^c$-structures are known to exist 
without any restriction on the topology of the underlying symplectic manifold 
(Forger, Hess \cite{FH} or Robinson, Rawnsley 
\cite{RR}). Compare approaches of 
Habermann, Habermann in \cite{HH} and  Cahen, Gutt, La Fuente-Gravy,
Rawnsley in \cite{CGLR}.

\bigskip

Similarly as in the orthogonal spinor case, 
we obtain by induction from representations $\sigma^k:\widetilde{G} \to 
\mbox{Aut}(H^k)$ and the principal bundle
$\mathcal{P} \to M$  bundles
$$p_k: \mathcal{H}^k = \mathcal{P} \times_{\sigma^k} H^k  \to M$$  where 
$k=0,\ldots, 2n.$ We call the bundle  $\mathcal{H}^{\bullet} = \bigoplus_{k=0}^{2n} 
\mathcal{H}^k$ the {\it higher oscillator bundle}.  
The basic oscillator bundle $\mathcal{H} = \mathcal{H}^0  =  \mathcal{P} 
\times_{\sigma} H$ coincides with the dual of the Kostant's or symplectic 
spinor bundle as it is called  by Habermann, Habermann \cite{HH}.
 Their continuity and 
differentiability are  discussed in the next sections.

\bigskip

{\bf Definition:}  The Banach bundle $\mathcal{CH}$ associated to $\mathcal{P}$ 
via  $\rho: 
\widetilde{G} \to \mbox{Aut}(CH),$ i.e., $$\mathcal{CH}  = 
\mathcal{P} \times_{\rho} CH$$ is called  the {\it bundle of measurements}. 

\bigskip 

The bundle $\mathcal{CH}$ exists regardless of the existence of the 
metaplectic structure because $\mbox{Ker}\, \lambda =\mbox{Ker}\, \rho.$

\bigskip

{\bf Remark} (Terminology related to measurements in Quantum theory):
The names `bundle of measuring devices' and `bundle of filters' for the bundle of 
measurements 
$\mathcal{CH}$ seem appropriate as well.  The  word
{\it filter} is used in quantum theory for a theoretical  
measuring device. See Ludwig \cite{Ludwig} for the role of 
filters, and still topical texts of
Birkhoff, von Neumann  \cite{BirkhoffNeumann} and  von Neumann \cite{NeumannKM}
for a mathematical framework for  quantum  (and also classical) 
 physical theories.

In each
point $m$ of a symplectic manifold $(M,\omega)$ (a `curved' phase space), we  
imagine 
a measuring 
device (filter) constituting of  filtering channels,  which are represented by 
elements of $\mathcal{CH}_m,$ the fiber in $m.$ 
Assuming that for any quantum measurement's result, 
there is a corresponding measuring device, $\mathcal{CH}_m$ presents a stack 
for all possible results of quantum measurements that can be 
obtained on a system found classically in the point $m.$ 
An ideal filter in point $m$ is a measuring device corresponding to a rank $1$ 
operator in $\mathcal{CH}_m$. Let $k: \mathcal{CH} \to M$ be the bundle 
projection of the bundle of measurements onto $M.$
Suppose that we have  measuring devices filling-up a submanifold 
$\iota: M' \hookrightarrow M$ (possibly only a discrete set of points).
 Then we have the pull-back bundle $k_{M'}: \iota^*(\mathcal{CH})  \to 
M'.$ In real measurements, each measuring device corresponds to a  
finite linear combination of rank 1
operators in the fibers $(\iota^*(\mathcal{CH}))_m, m\in M'.$

\subsection{Continuity of $G$-structures}

We make a comment on a canonical continuity of atlases of induced Banach 
bundles.
Let us recall that for a  principal 
$Q$-bundle $\pi: \mathcal{R} \to M$
and a  representation $\kappa: Q \to \mbox{Aut}(X)$ of 
a Lie group $Q$ on a topological vector space $X$, the induced 
(associated) bundle is the quotient $\mathcal{F}=(\mathcal{R} \times X) 
/\simeq$ together with the appropriate bundle projection onto $M$. 
The equivalence 
relation $\simeq$ is defined 
for any $r, s  \in \mathcal{R}$ and $v,w \in X$ by
$(r,v) \simeq (s,w)$ if and only if there exists an element $g\in Q$ such that
$s = r g^{-1}$ and $w=\kappa(g)v.$ The topology induced on the quotient  
is the final topology, i.e., 
the finest one for which the canonical projection $q: \mathcal{R} \times X \to 
(\mathcal{R} \times X)/ \simeq$ is continuous.
The bundle projection $p:\mathcal{F} \to M$ is given by $p([(r,v)])=\pi(r),$ 
where $[(r,v)] \in \mathcal{F}.$
For $U\subseteq M,$ let $\mathcal{F}_U$  denote the set $p^{-1}(U).$ 
Especially, 
$\mathcal{F}_{\{m\}},$ denoted by $\mathcal{F}_m,$ is the fiber of 
$\mathcal{F}$ 
in $m\in M.$

Let us equip the higher oscillator bundle and the bundle of measurements 
with the final topologies and prove that they form Hilbert and Banach 
bundles, respectively. Actually, it is convenient to  prove a more general 
theorem.  

\bigskip

{\bf Theorem 4:} Let $\kappa: Q \to \mbox{Aut}(X)$ be an isometric or a unitary
representation of a Lie group $Q$  on Banach or Hilbert space $X,$ 
respectively.  Let $p:\mathcal{F} \to M$ be the associated bundle to a 
principal $Q$-bundle $\pi: \mathcal{R} \to M.$ Then $p$ is canonically a 
continuous Banach or a Hilbert bundle, respectively.
 
{\it Proof.} Let us take a principal bundle atlas 
$(U_{\alpha},\phi_{\alpha})_{\alpha \in I}$ in the differentiable structure of 
the
principal bundle $\pi.$  In particular, $\phi_{\alpha}: \pi^{-1}(U_{\alpha}) 
\to U_{\alpha} \times Q.$
For to get a bundle atlas on $p:\mathcal{F} \to M,$ we define  maps 
$\psi_{\alpha}': U_{\alpha} \times X \to p^{-1}(U_{\alpha})$ 
by $\psi_{\alpha}'(m,v)= q(\phi^{-1}_{\alpha}(m,e),v),$ where $e$ is the unit 
element in $Q,$ $m \in U_{\alpha}$ and $v \in X.$ The maps $\psi'_{\alpha}$ 
are 
continuous because
$\phi^{-1}_{\alpha}$ and $q$ are continuous. They are fiber preserving in the 
sense that $p \circ \psi_{\alpha}' = \mbox{pr}_1,$ where 
$\mbox{pr}_1: U_{\alpha} \times X \to U_{\alpha}$ is the projection onto the
first factor in the product. Their inverse maps, $\psi_{\alpha}: 
p^{-1}(U_{\alpha}) \to
U_{\alpha} \times X$ are given by $\psi_{\alpha}([(r,v)])=(\pi(r), 
\kappa(g^{-1})v)$ where $g$ is the unique element of $Q$ such that 
$rg=\phi_{\alpha}^{-1}(m,e)$ 
and  $m=\pi(r).$ The independence on the chosen representative of an 
equivalence 
class 
follows from the definition of the associated bundle. Let $U \subseteq 
U_{\alpha}$ and $Y \subseteq X$ be open sets.
The $\psi_{\alpha}$-preimage of $U \times Y$ is 
$S=\{[(\phi^{-1}_{\alpha}(m,e),w)]|\, m \in U,  w \in Y\}.$ 
Consequently, the $q$-preimage $q^{-1}(S) = \{(\phi_{\alpha}^{-1}(m,e)g^{-1}, 
\kappa(g)w)|\, m \in U, g \in Q,   w \in Y\}.$
The preimage can be written as $\bigcup_{g\in Q}(\phi_{\alpha}^{-1}(U,e)g^{-1}, 
\kappa(g)Y).$\footnote{More properly, the image equals to
$\bigcup_{g \in Q}\{(\phi_{\alpha}^{-1}(m,e)g^{-1}, 
\kappa(g)y)|\, m \in U, y \in Y\}.$} 
For each $g\in Q,$ both of the sets $\{\phi_{\alpha}^{-1}(m,e)g^{-1}|\, m \in 
U \}$ and $\{ \kappa(g)y|\, y \in Y \}$ are open since  $\phi_{\alpha}^{-1}$ 
is a homeomorphism and   the map $\kappa(g)$ is a homeomorphism  as 
well. Due to the definition of the topology on Cartesian products, $q^{-1}(S)$ 
is 
open and thus $\psi_{\alpha}$ is continuous.
Summing-up, $\psi_{\alpha}$ is a homeomorphism and  
$\mathfrak{A}=(U_{\alpha},\psi_{\alpha})_{\alpha \in I}$ is  
a bundle atlas.

Choosing a different atlas in the differentiable structure of 
$\mathcal{R},$ 
leads to an induced atlas $\widetilde{\mathfrak{A}}=(\widetilde{U}_{\alpha}, 
\widetilde{\psi}_{\alpha})_{\alpha \in \widetilde{I}}$ that is  different from  
$\mathfrak{A}$ in general. However, the transition functions $\psi_{\alpha 
\beta}$ for $\mathfrak{A} 
\cup \widetilde{\mathfrak{A}}$ are given by 
composition of the representation with the 
transition functions $(\phi_{\alpha \beta})_{\alpha, \beta}$ of the principal 
bundle, i.e., $\psi_{\alpha \beta} = \kappa \circ \phi_{\alpha \beta}.$
Since $\kappa$ is continuous (as a map $Q \times X \to X)$ and unitary, it is 
continuous
as a map $Q \to U(X)$ where $U(X)$ is considered with the strong topology (see 
Schottenloher \cite{Schotten} or
 Knapp \cite{KnappOverview}, pp. 10 and 11). Thus we have defined an independent continuous 
bundle structure.

Suppose that $\kappa$ is a unitary representation on a Hilbert space $X$. The 
fiber-wise scalar product is defined by $([(r,v)],[(r,w)])_S = (v,w)_X$ 
where at the right-hand side, the scalar product on $X$ is meant. The 
independence of the scalar product on the chosen representatives follows from 
the fact that the representation is unitary.
Let $|| \,||_S:\mathcal{F} \to \mathbb{R},$ given by $||a||_S=\sqrt{(a,a)_S},$ 
$a\in \mathcal{F},$ be the induced norm, and let $(a_n)_n$ converge to an 
element $a$ in the 
quotient topology of $\mathcal{F}.$  By definition of the quotient topology, 
there exists a 
selector $(b_n)_n$ on $\{q^{-1}(\{a_n\})|\, n \in \mathbb{N}\}$ which converges 
in 
$\mathcal{R} \times X.$ Denoting its limit by $b,$ we have that 
1) $c_n=\mbox{pr}_2(b_n)$ converges to $c=\mbox{pr}_2(b)$  in $X$, and 2)
$a=[b],$ where $\mbox{pr}_2$ denotes the projection $\mathcal{R} 
\times X \to X$ onto the second factor. Summing-up, 
$||a_n||_S=\sqrt{(a_n,a_n)_S}=\sqrt{(c_n,c_n)_X} 
\rightarrow ||c||_X=||[b]||_S =||a||_S.$ Thus, $||\,||_S$ is continuous.
The continuity of $(,)_S$ in each of the arguments follows from the 
polarization 
identity.

Rewriting the definition of $(,)_S$, we have
$(q(r,v),q(r,w))_S=(v,w)_X.$
Further $(\psi_{\alpha}'(m,v),\psi_{\alpha}'(m,w))_S=
(q(\phi_{\alpha}^{-1}(m,e),v), q(\phi_{\alpha}^{-1}(m,e),w))_S=$ 
$(v,w)_X$, \\ proving the unitarity of $\psi_{\alpha}'.$ 
Consequently,  $\mathcal{F}$ is a Hilbert bundle. 

Similarly one proceeds in the case of the isometric 
Banach representation. 

\hfill\(\Box\)

\bigskip

Let us call any atlas $(U_{\alpha}, \psi_{\alpha})_{\alpha \in I}$ on an 
induced 
bundle $\mathcal{F}$ constructed as in the preceding proof 
the {\it canonically induced  atlas} or, in  more detail, the 
atlas canonically induced 
by the principal bundle atlas $(U_{\alpha}, \phi_{\alpha})_{\alpha \in I}.$ As 
follows from the above proof,
canonical atlases induced by different but compatible 
atlases 
are $C^0$-compatible. In particular, one can speak about the induced 
$C^0$-differentiable structure.

\bigskip

{\bf Corollary 5:} For a symplectic manifold $(M,\omega)$ and a  metaplectic 
structure
$(\pi_P, \Lambda)$ on it, the higher oscillator bundle  
and the bundle
of  measurements, both equipped with the canonically induced atlases by a 
principal bundle atlas, are continuous Banach bundles. The basic oscillator 
bundle is a 
continuous Hilbert bundle.

\subsection{Differentiable structures on oscillator bundle}
  
A substantial obstacle for analysis on induced Banach or Hilbert bundles of
infinite rank is that they need  not be smooth or not even   $l$-times
 differentiable for some $l \in \mathbb{N}.$
When speaking about $C^l$-sections, one has to take care to which atlas the 
degree of differentiability refers.

\bigskip

Let us recall that for a representation $\kappa: Q \to \mbox{Aut}(X)$ of a Lie 
group $Q$ on a 
Banach space $X,$ a vector $v \in X$ is called smooth if
$Q\ni g \mapsto \kappa(g)v \in X$ is smooth, i.e., it has each Fr\'echet 
differential in each point $g.$

\bigskip

Let  $\phi_{\alpha \beta}: U_{\alpha} \cap U_{\beta} \to  Q$ be a 
cocycle of  transition functions of a principal $Q$-bundle $\mathcal{R}$ 
corresponding to a
smooth atlas $(U_{\alpha}, \phi_{\alpha})_{\alpha \in I}$ of $\mathcal{R}.$ 
Using  the definition of the canonically induced atlas, we obtain that its
transition maps $\psi_{\alpha \beta}: (U_{\alpha} \cap U_{\beta}) \times X \to 
(U_{\alpha} \cap U_{\beta}) \times X$ are given by
$\psi_{\alpha \beta}(m,v) = (m, \kappa(\phi_{\alpha \beta}(m))v).$
We agree that a Banach bundle atlas is called smooth if the transition 
functions  $\psi_{\alpha \beta}$ considered as maps of $U_{\alpha} \cap 
U_{\beta} \to \mbox{Aut}(X) \subseteq \mbox{End}(X)$ are 
smooth. On $\mbox{Aut}(X),$ we consider the topology 
induced from the strong operator topology on $\mbox{End}(X).$
This happens if and only if $m \mapsto 
\kappa(\phi_{\alpha \beta}(m))$ is 
smooth. For that it is sufficient the representation to be smooth as a map into
$\mbox{End}\,(X)$ with the strong topology.

Thus we get

\bigskip

{\bf Lemma 6:} The  bundle atlas canonically induced by a principal bundle 
atlas 
is smooth if the inducing 
representation is 
smooth.

\bigskip

The opposite implication does not hold. For it, take an arbitrary 
manifold $M$, and the metaplectic group $Mp(V,\omega)$ for the Lie group $Q.$ 
Let us consider the product bundle $M \times Q \to 
M$ as the principal $Q$-bundle and  the oscillator representation $\sigma$ 
as the representation $\kappa.$
For the smooth atlas on the product bundle take the one containing the 
single chart $U = M$ and the only bundle map $\phi(m,g)=(m,g),$ where $m\in M$ 
and $g\in Q.$ The set of transition functions of this atlas consists just of
the constant map $m \mapsto e,$ where $e$ is the unit element in the metaplectic group. 
The set of transition functions of the canonically induced atlas on 
$\mathcal{F} = 
(M \times Q) \times_{\sigma} \check{H}$ contains only $\psi(m,v)=(m,\kappa(e)v) 
= (m,v)$, which is clearly smooth, although the 
representation is known not to be. (See Borel, Wallach  \cite{BorelWallach}, p. 
153, or below.)

\bigskip

{\bf Smooth Hilbert bundle \`a la Kuiper} 

\bigskip

  Because the unitary group of $H$ is contractible in the strong operator
topology (Dixmier, Douday \cite{DixmierDouady})\footnote{It is contractible 
also in the operator norm topology (Kuiper \cite{Kuiper}).},  
  there exists an isomorphism  $J: \mathcal{H} \to 
M  \times H$ of topological vector bundles.  We declare the smooth structure 
of $\mathcal{H}$ as the 
maximal set of smoothly compatible atlases containing the chart  
$(M,J)$ and 
call it the {\it Kuiper structure} induced by $J$.

\bigskip

{\bf Smooth Fr\'echet bundle \`a la Bruhat}

\bigskip

Let $S(L)$ denote the Schwartz space of rapidly decreasing smooth functions on 
$L.$ Vectors in $\check{H} \setminus S(L)$ are exactly the ones which are 
not smooth for $\check{\sigma}$ (Borel, Wallach \cite{BorelWallach}).  
 Consequently, smooth vectors of $H^k$ acted  on by $\sigma^k$ form
 spaces $S^k = \bigwedge^{k}V^* \otimes 
S,$ where $S \subseteq H$ denotes the $\flat$-image of $S(L) \subseteq 
L^2(L)$ in $L^2(L)'.$
Note that
$S$ is not complete in the norm topology inherited from $H.$  To get a complete space, we 
retopologize it as it is usual in representation theory (see below).
  
We recall and combine some known facts from Lie groups' representations (based on 
the  Bruhat's  dissertation thesis \cite{Bruhat})
and some known facts specific for the Segal--Shale--Weil representation.

\begin{itemize}
\item[1)] Representations $\sigma^k$ preserve the space of smooth vectors 
$S^k.$
\item[2)] Let us denote the restrictions  $\sigma^k$ to $S^k$ by $\sigma^k_0$ 
and consider $S^k$ to be equipped with the topology inherited by  the
injection $S^k \to C^{\infty}(\widetilde{G}, H^k),$ $S^k \ni v \mapsto (\pi^v:\widetilde{G} \ni g \mapsto \sigma^k(g)(v) \in H^k).$
Here, $C^{\infty}(\widetilde{G}, H^k)$ is considered with the compact $C^{\infty}$-topology. 
(See,  e.g., Warner \cite{Warner}, p. 484 for a definition of 
compact $C^{\infty}$-topology if necessary.)  The restricted representation 
$\sigma_0^k$ is continuous in the new topology (\cite{Warner}, Sect. 4.4.). 

\item[3)] The compact $C^{\infty}$-topology on $S$ coincides with the  classical
semi-norm Fr\'echet topology on this space (see Borel, Wallach \cite{BorelWallach}). 

\item[4)] Since (a) $\sigma^k_0$ is a  representation on $S^k$ with 
the compact $C^{\infty}$-topology, which is Fr\'echet, and since (b) the space 
$(E_{\infty})_{\infty}$ of smooth vectors of the space $E_{\infty}$ of smooth 
vectors of a continuous representation 
on a Fr\'echet space $E$\footnote{Both $E_{\infty}$ and $(E_{\infty})_{\infty}$ 
are considered with the compact $C^{\infty}$-topology.} is topologically 
isomorphic to $E_{\infty},$ we conclude that $\sigma_0^k: \widetilde{G} \to 
\mbox{Aut}\,(S^k)$ are smooth maps
when $\mbox{Aut}(S^k) \subseteq \mbox{End}(S^k)$ is considered with the strong 
topology. 
Technically speaking, this is because the assumption of Thm. 4.4.1.7 
(see Warner \cite{Warner}) is satisfied, and thus Remark (3) at pages 258 and 259 
of \cite{Warner} applies.

\end{itemize}

Taking representations $\sigma_0^k$ instead of $\sigma^k,$  we get
Fr\'echet bundles $\mathcal{S}^k= \mathcal{P} \times_{\sigma^k_0} S^k.$  
By  what was mentioned in item 4 above and by Lemma 6, each of the 
  bundles $\mathcal{S}^k$ are smooth with respect to the atlas induced 
canonically by any smooth principal bundle atlas (on $\mathcal{P}).$
We call the bundles $\mathcal{S}^k \to M$ with the smooth atlases {\it Bruhat structures}. The atlases
are constructed as in the proof of Theorem 4.
 
\bigskip

{\bf $C^0$-Hilbert bundle \`a la Goodman}

\bigskip

We can consider  spaces $W^l$ of $C^l$-vectors in $H$ as well, $l \geq 
1,$ obtaining  representations $\sigma^l: \widetilde{G} \to  
\mbox{Aut}(W^l)$.  See Goodman \cite{Goodman} or Neeb 
\cite{Neeb}. Consider the tensor product representations 
$\sigma^{k,l}:\widetilde{G} \to \mbox{Aut}(\bigwedge^k V^* \otimes W^l).$  

\bigskip

 From  the  known form of the differential of the Segal--Shale--Weil
representation (see Habermann \cite{HabermannAn}, or Kirillov \cite{Kirillov}, 
p. 184 for the case $n=1$) and from  the definition of the so-called 
Hermite--Sobolev spaces (as defined, e.g., in Bongioanni, Torrea \cite{Bong}), 
spaces $W^l$ are isomorphic to the  Hermite--Sobolev spaces denoted 
by $W^{2l}_2$ or sometimes $W^{2l,2}$ in \cite{Bong}.
 Moreover, the compact $C^l$-topology on $W^l,$ as described in 
Goodman \cite{Goodman}, coincides with the norm topology, generated by the oscillator ladder operators, considered in Bongioanni, Torrea \cite{Bong}.
(Notice that it is easy to see that spaces $W^l,$ with this topology, are 
Hilbert spaces as well as that  $\sigma^{0,l}$ is not $C^1.$)

\bigskip

The associated bundles $\mathcal{W}^{k,l} =\mathcal{P} 
\times_{\sigma^{k,l}}(\bigwedge^k V^* \otimes W^l)$ when equipped with the 
canonically induced atlases are Hilbert 
bundles with in general, only continuous structures.  We call
these bundles $\mathcal{W}^{k,l} \to M$ with the continuous atlases the 
{\it Goodman structures}.

\bigskip

{\bf Remark:} In the case of the Bruhat or the Goodman structures, we 
 either loose the completeness with 
respect to a norm or the smoothness of the representation, respectively.

\subsection{Differentiable structure on bundle of measurements}

Operations $\sharp: H \to \check{H}$ and $\flat: \check{H} \to H$ lift to  
bundles
$\mathcal{H}\to M$ and $\check{\mathcal{H}} \to M,$ respectively. We denote 
them by the same symbols, i.e.,
$\sharp: \mathcal{H} \to \check{\mathcal{H}}$ and $\flat: \check{\mathcal{H}} 
\to \mathcal{H}.$
We shall denote the operator $\mbox{Id}_M \times \sharp: M \times H \to M 
\times 
\check{H}$ by $\sharp_M$
and the operator $\mbox{Id}_M \times \flat: M \times \check{H} \to M \times H$ 
by $\flat_M.$
 Let us consider the map
$$\otimes_{\mathcal{CH}}: \check{\mathcal{H}} \times_{M} \mathcal{H} \to 
\mathcal{CH} \mbox{ defined by } \, [(r,v)]\otimes_{\mathcal{CH}} [(r,f)] = 
[(r, v \otimes f)]$$ where $v\in \check{H}, f\in H$ and $r\in \mathcal{P}.$  
Since the operator norm topology on $CH$ is coarser than
 the Hilbert tensor product topology on $\check{H} \otimes H,$  described 
equivalently by the Hilbert--Schmidt norm,  
$\otimes_{\mathcal{CH}}$ is 
continuous. For a moment, let us denote by $\underline{\mathbb{C}}$ the trivial 
complex line 
bundle $M \times \mathbb{C}$ over $M$. This bundle is  isomorphic to 
the
bundle associated to $\mathcal{P}$ by the 
trivial representation of $\widetilde{G}$ on $\mathbb{C}.$ Further, we set
\begin{align*}
ev_{\mathcal{H}}: \mathcal{H} \times_M \check{\mathcal{H}} \to 
\underline{\mathbb{C}}& \qquad
ev_{\mathcal{H}}([(r,f)],[(r,v)])= (\pi_P(r),f(v))\\
\end{align*}
Instead of  
$\mbox{pr}_2\left(ev_{\mathcal{H}}([(r,f)],[(r,v)])\right),$ we write 
$[(r,f)]([(r,v)]),$ where $\mbox{pr}_2: M\times 
\mathbb{C} \to \mathbb{C}$ denotes the projection onto the second factor.
The operations $\sharp: \mathcal{H} \to \check{\mathcal{H}}, \flat:\mathcal{H}
\to \check{\mathcal{H}}, \otimes_{\mathcal{CH}}:\check{\mathcal{H}}\times_M 
\mathcal{H} 
\to \mathcal{CH}$ and  $ev_{\mathcal{H}}:\mathcal{H}\times_M\check{\mathcal{H}} 
\to\underline{\mathbb{C}}$ are 
well defined. The proof of this is based on the $\widetilde{G}$-equivariance 
(Lemma 3) and goes in the same lines as the one of the 
correctness of definitions of the $\mathcal{CH}$-product and the
$\mathcal{CH}$-action, explained below (Sect. 4. 1).

\bigskip

{\bf Notation for bundle $\underline{H}:$} For the trivial bundle 
$\underline{H} = M \times H \to M,$ we set $(m,f')+_M (m,f'') = (m,f'+f''),$  
$t\odot_M (m,f) = (m,tf),$ and $((m,f'),(m,f'')) = 
(f',f'')_H \in \mathbb{C}$
where $t\in \mathbb{R},$ $f,f',f'' \in H,$  and $m\in M.$
This notation may look like superfluous. However, we use 
it later (sect. 4.3) when it appears to be helpful.

\bigskip

{\bf Normalized Kuiper maps:} 
Since $\sigma$ is a unitary representation of $\widetilde{G}$ on $H,$ the 
associated bundle $\mathcal{H}$ is a vector bundle with fiber $H$ and structure 
group
$U(H)$ with the strong operator topology (Thm. 4).
The set of isomorphism classes of vector bundles with fiber $H$ and structure group $U(H)$ is in a set bijection
with the first \v{C}ech cohomology group $H^1(M, U(H))$ (see, e.g., Raeburn, Williams \cite{Raeburn}, pp. 96--97).
Because the unitary group of an infinite dimensional Hilbert space is  contractible,
the \v{C}ech cohomology group consists only of one element. Consequently, each bundle must be isomorphic to
the trivial one. This means that by definition (p. 93, \cite{Raeburn}), there 
is a fiber preserving homeomorphism $J_N$ of $\mathcal{H}$ and $\underline{H}$
which is a map into $U(H)$ when considered composed with  bundle maps. In particular, this  implies that $J_N$ is a unitary map in each fiber.
We call such a map a  normalized Kuiper map. {\bf By a Kuiper map, we shall only mean a normalized Kuiper map}, i.e., 
a unitary isomorphism.

\bigskip

Since  the bundle of measurements $\mathcal{CH}$ is a specific associated 
bundle to the  principal bundle $\mathcal{P},$  
the transition functions of $\mathcal{CH}$ are compositions of 
the transition functions of  $\check{\mathcal{H}},$ the representation 
  $\sigma,$ and the $\mbox{Ad}$-representation. Consequently, $\mathcal{CH}$ is a 
trivial bundle  by  the Dixmier--Douady theory. (See Raeburn, Williams \cite{Raeburn},  Thm. 4.85 
(a), pp. 109--110.)  
Because the cited result is often scattered through texts,
we give an elementary proof based only on the triviality of 
$\mathcal{H}$ and $\check{\mathcal{H}}$ and  
functional analysis' methods.

\bigskip

{\bf  Theorem 7:} Let $(M,\omega)$ be a symplectic manifold admitting a 
metaplectic structure $(\pi_P, \Lambda).$
The bundle of measurements  $\mathcal{CH} \to M$ and the 
product bundle $M \times CH \to M$ are isomorphic as Banach bundles. 
Consequently, the bundle of measurements is trivial.

{\it Proof.} 
 Let us consider the  Banach bundle $\check{\mathcal{H}}=
\mathcal{P} \times_{\check{\sigma}} \check{H}$ and let $J$ be a Kuiper map of 
the bundle $\mathcal{H} \to 
M$, i.e., $J: \mathcal{H} \to M \times H$ is an isomorphism of continuous 
Hilbert bundles which is isometric on each fiber.
 It induces a map $\check{J}: \check{\mathcal{H}} \to M \times \check{H}$ 
(trivializing $\check{\mathcal{H}}$)
given by the formula
$\check{J}(\tilde{v}) = (J(\tilde{v}^{\flat}))^{\sharp_M},$ where $\tilde{v} 
\in 
\check{\mathcal{H}}.$ 
Because $J$, $\sharp_M$ and $\flat$ are continuous, so is $\check{J}.$
The inverse to $\check{J}$ is given by $\check{J}^{-1}(v) = 
(J^{-1}(v^{\flat_M}))^{\sharp},$ $v \in M \times \check{H}.$
Since $\flat_M, J$ and $\sharp$ are continuous, we get that $\check{J}$ 
is a homeomorphism and thus $(M,\check{J})$ is a bundle chart for 
$\check{\mathcal{H}}.$

 For bundle rank 1 operators, we set 
$$J_{CH}^0([(r, v \otimes f)])= 
\left(\pi_P(r),\check{\mbox{pr}}_2(\check{J}([(r,v)])) \otimes 
\mbox{pr}_2(J([(r,f)]))\right)$$ 
where $f \in H,$ $v \in \check{H},$
$r \in \mathcal{P},$ $\mbox{pr}_2: M \times H \to H$ and $\check{\mbox{pr}}_2: 
M 
\times  \check{H} \to \check{H}$ denote the
projections onto the second factors of the respective products.  
To prove that  $J_{CH}^0$ is well defined, we use formula (\ref{formule}) from 
Lemma 3. Indeed, 
\begin{align*}
J_{CH}^0([(rg ^{-1}, \rho(g)&(v \otimes f))])= J_{CH}\left([(rg^{-1}, 
\check{\sigma}(g)v \otimes \sigma(g)f)]\right)\\
& = \left(\pi_P(rg^{-1}), 
\check{\mbox{pr}}_2\check{J}([(rg^{-1},\check{\sigma}(g)v)]) \otimes 
\mbox{pr}_2J([(rg^{-1},\sigma(g)f)])\right)\\
& =\left(\pi_P(r), \check{\mbox{pr}}_2 \check{J}([(r,v)]) \otimes 
\mbox{pr}_2J([(r,f)])\right)  = J_{CH}^0([(r, v \otimes f)])
\end{align*}
Since  defined by composition of continuous maps, $J_{CH}^0$ is continuous.

The inverse of $J_{CH}^0$ on rank 1 operators is given by
$(J_{CH}^0)^{-1}(m,\tilde{v} \otimes \tilde{f}) = \check{J}^{-1}(m,\tilde{v}) 
\otimes_{\mathcal{CH}} J^{-1}(m,\tilde{f}),$ $m \in M,$ $\tilde{v} \in 
\check{H}$ and 
$\tilde{f} \in H,$   which is a composition of continuous maps.
On finite rank operators, $J^0_{CH}$ is extended linearly. Obviously, 
$(J_{CH}^0)^{-1}$ on finite rank operators is 
continuous as well.

Finally for $r\in \mathcal{P}$ and $a \in CH,$ we set 
$$J_{CH}([(r,a)])=\left(\pi_P(r),\lim_i \mbox{pr}_2^{CH}J_{CH}^0([(r, 
a_i)])\right)$$ where 
$(a_i)_{i\in \mathbb{N}}$ is any 
sequence of finite rank operators converging
to $a$ in the operator norm topology, and $\mbox{pr}_2^{CH}$ denotes the 
projection of 
$M 
\times CH$ onto the second factor.  The map $J_{CH}$ is independent of
the chosen sequence. Indeed, let $(a_i)_{i\in \mathbb{N}}$ and $(b_i)_{i\in 
\mathbb{N}}$ 
be sequences of finite rank operators which converge to $a$ in the operator 
norm. 
Their difference $c_i = a_i - b_i$ converges to the zero operator.
Moreover, each $c_i$ is a finite rank operator. Therefore, there are 
$\lambda^{jk}_i \in \mathbb{C}$ such that
$c_i=\sum_{j,k=1}^{n_i}\lambda^{jk}_i h_j \otimes \eta_k,$ where $(h_j)_j$ and 
$(\eta_k)_k$ are complete orthonormal systems of $\check{H}$ and $H,$
respectively. We have $||c_i||_{CH}^2 = \sum_{j,k}|\lambda^{jk}_i|^2.$ Since 
$c_i 
\to 0,$ we get $\lambda^{jk}_i \to 0$ for all $j,k.$
This fact and the linearity of $J$ and $\check{J}$ make us able to conclude 
that $J_{CH}^0([(r,c_i)])=\lim_i 
\sum_{j,k}\lambda^{jk}_i\left(\check{\mbox{pr}}_2\check{J}([(r,h_j)]) \otimes 
\mbox{pr}_2 J([(r,\eta_k)])\right)=
0.$ This proves that the definition of $J_{CH}$ is correct.

Map $J_{CH}$ is continuous since it is defined as the continuous extension in 
the 
second argument and 
since $\pi$ is continuous.
The inverse of $J_{CH}$ is given by $J_{CH}^{-1}(m,a)=\lim_i 
(J^{CH}_0)^{-1}(m,a_i),$ 
where $(a_i)_{i\in \mathbb{N}}$ is a sequences of finite rank operators 
converging 
to $a$ in the operator norm.
Summing-up, $J_{CH}: \mathcal{CH} \to M \times CH$ is a homeomorphism of 
 $\mathcal{CH}$ and  $M\times CH.$ Since $J_{CH}$ is fiber preserving, $J_{CH}$ 
is an 
isomorphism of continuous bundles. Since the Kuiper maps are taken to be 
normalized, $J_{CH}$
is a fiber-wise isometry and thus, it is an isomorphism of Banach bundles and  
consequently, $\mathcal{CH}$ is trivial as a Banach bundle.
\hfill\(\Box\)

\bigskip

{\bf Remark:}  It is known that the set  of equivalence classes of  
bundles of Banach spaces, with  fibers the  algebra of 
compact operators $CH$ on a Hilbert space $H$ and 
whose structure group is the projective unitary group of $H,$  is isomorphic to 
the third 
singular cohomology 
group $H^3(M,\mathbb{Z})$ 
of $M.$ The equivalence is the isomorphism of continuous fiber bundles.
 See  Dixmier, Douady \cite{DixmierDouady} or Dixmier \cite{Dixmier}. 
We refer also to Mathai, Melrose and Singer \cite{MMS} and Schottenloher 
\cite{Schotten} for a context and for a nice topological study, 
respectively. (Moreover, the set of such Banach fiber bundles 
forms a group under the spatial tensor product (the Brauer group). This group 
is isomorphic to the additive  structure on $H^3(M,\mathbb{Z}).$ See 
Rosenberg \cite{Rosenberg}.)

\section{Hilbert $C^*$-bundle  and  Kuiper complex} 

In this section, we first introduce pointwise structures which we use for 
defining  a Hilbert $C^*$-structure on the higher oscillator bundle. 
Then we use a product type connection to establish the  
elliptic complex, mentioned in the Introduction.

\subsection{Pointwise analytical structures}
 
 Let us define a  map  $\cdot_{\mathcal{H}}:  
\mathcal{H}^{\bullet} \times_M \mathcal{CH} \to 
\mathcal{H}^{\bullet}$ by 
$$[(q,f)] \cdot_{\mathcal{H}} [(q,a)] = [(q, f \cdot a)]  $$ 
where $q \in \mathcal{P},$ $f\in 
H^{\bullet},$ and  $a\in CH.$ We call it the {\it $\mathcal{CH}$-action}.
The independence on the representative of an equivalence class is proved  below.

\bigskip

For to define a pointwise $\mathcal{CH}$-valued product on 
$\mathcal{H}^{\bullet},$ we 
reduce the structure group of $\mathcal{P}$ 
from $\widetilde{G}$ to $\widetilde{K}$ which is possible since  this group
is  a deformation retract of the structure group $Mp(V,\omega_0)$ by the 
lifting 
property for coverings
 and the known fact that $Sp(V,\omega_0)$ deformation retracts onto
$K$.   
The reduction is non-unique. Therefore, we choose a compatible positive almost 
complex structure on $(M,\omega)$ and consider the $\Lambda$-preimage of the bundle of 
unitary 
frames. We denote the resulting principal 
$\widetilde{K}$-bundle by $\mathcal{P}_R.$  
The corresponding 
associated bundles are denoted by $\mathcal{H}^{\bullet}_{R}$ and 
$\mathcal{CH}_{R}.$ We call the former one the  reduced higher oscillator 
bundle. The total spaces of these vector bundles  
coincide,  as topological spaces, with $\mathcal{H}^{\bullet}$ and 
$\mathcal{CH},$ respectively. 
Using this reduction, we  define the pointwise {\it $\mathcal{CH}$-product}
$(,)_{\mathcal{CH}}: \mathcal{H}^{\bullet}_R \times_M \mathcal{H}^{\bullet}_R 
\to 
\mathcal{CH}_R \mbox{  by the  formula }$
$$([(q,v)],[(q,w)])_{\mathcal{CH}}=[(q,(v,w))] \in \mathcal{CH}_R$$ where $q 
\in 
\mathcal{P}_R,$  and $v,w \in H^{\bullet}.$ 
Note that whereas at the right-hand side, $(v,w) \in CH,$ at the left-hand 
side, a map
$(,)_{\mathcal{CH}}:\mathcal{H}^{\bullet}_R \times_M \mathcal{H}^{\bullet}_R 
\to 
\mathcal{CH}_R$ is prescribed as we prove yet.

\bigskip 
 
{\bf Lemma 8}: For a symplectic manifold $(M,\omega)$ admitting a metaplectic 
structure $(\pi_P,\Lambda)$,  the maps
\begin{eqnarray*}
\cdot_{\mathcal{H}}&:&  \mathcal{H}^{\bullet} \times_M \mathcal{CH} \to 
\mathcal{H}^{\bullet}\\
(,)_{\mathcal{CH}}&:& \mathcal{H}^{\bullet}_R \times_M \mathcal{H}^{\bullet}_R 
\to 
\mathcal{CH}_R
\end{eqnarray*}
 are  correctly defined.

{\it Proof.} 
We have to check that these  maps do not depend on the choice of 
representatives.
\begin{itemize}
 
\item[1)]   Let us verify that the definition of the action 
$\cdot_{\mathcal{H}}$ 
is 
correct. For $q\in \mathcal{P},$ $g \in 
\widetilde{G},$ $k=0,\ldots, 2n,$  $a\in CH,$ and 
$v = \alpha \otimes h \in H^{\bullet}$ 
\begin{align*}
[(q,v)]\cdot_{\mathcal{H}} [(q g, \rho(g)^{-1}a)]  &= [\left(q g, 
\sigma^k(g^{-1})v\right)]\cdot_{\mathcal{H}} [\left(q 
g,(\rho(g)^{-1}a)\right)]\\
&= [\left(q g, \left(\sigma^k(g)^{-1}(v)\right) \cdot  
\left(\rho(g)^{-1}a\right)\right]\\
&= [\left(q g, \sigma^k(g)^{-1} (v \cdot a)\right)]\\
&=  [(q,v)]\cdot_{\mathcal{H}}[(q,a)] 
\end{align*}
where the equivariance of the $CH$-action (Lemma 3) is used in the second row.
In a similar way, one proves that the definition of $\cdot_{\mathcal{H}}$ does 
not depend on the choice of a 
representative of the element acted upon.

\item[2)] We prove  the correctness of the definition of $(,)_{\mathcal{CH}}.$ 
For a frame $q \in \mathcal{P}_R,$ $\alpha, \beta \in \bigwedge^{k}V^*,$ $v,w 
\in H$ and $g\in \widetilde{K}$  
\begin{align*}
\left([(q g,\sigma^k(g^{-1})v)],[(q , 
w)]\right)_{\mathcal{CH}}&=\left([\left(qg,\sigma^k(g^{-1})v\right)], [\left(q 
g,\sigma^k(g^{-1})w)\right]\right)_{\mathcal{CH}}\\
&=\left[\left(qg, \left(\sigma^k(g^{-1})v, 
\sigma^k(g^{-1})w\right)\right)\right]\\
&=\left[\left(q g,  \rho(g)^{-1}(v,w)\right)\right] = [\left(q,(v,w)\right)]
\end{align*}
In the second row above, the $\widetilde{K}$-equivariance of the 
$CH$-product from Lemma 3 is used.
The independence on representatives of the  second argument of 
$(,)_{\mathcal{CH}}$   follows from 
the hermitian symmetry of $(,): H^{\bullet} \times H^{\bullet} \to CH$.
\end{itemize}
\hfill\(\Box\)

\bigskip

{\bf Remark:} 
Consider the composition

$$\xymatrix{ 
h \in \mathcal{H}    \ar@{|->}[r]^-{(,)_{\mathcal{CH}}} & 
(h,h)_{\mathcal{CH}} \in \mathcal{CH} \ar@{|->}[r]^-{||\,||_{\mathcal{CH}}} &  
||(h,h)_{\mathcal{CH}}||_{\mathcal{CH}} \in \mathbb{R}_{\geq 0}
}$$

where $||[(q,a)]||_{\mathcal{CH}} = ||a||_{CH}$ for $q\in \mathcal{P}_R$ and 
$a\in 
CH.$
This map is well defined since $\rho$ is an isometry. Let us recall that there 
is  also a fiber-wise scalar product 
$(,)_S$ on $\mathcal{H}\to M$ as constructed in the proof of Thm. 4
that makes  $\mathcal{H}$ a Hilbert bundle.    One can prove easily 
that 
$||h||_S^2=||(h,h)_{\mathcal{CH}}||_{\mathcal{CH}}$ where $||\,||_S$ is the 
fiber-wise norm induced by $(,)_S.$

\subsection{Hilbert $C^*$-bundles}
 
According to Fomenko, Mishchenko \cite{FM}, we shall define a 
smooth action of  $CH$ on the total space $\mathcal{H}^{\bullet}_R,$ which 
restricts to a Hilbert 
$CH$-module action on each fiber; and a smooth
$CH$-valued map  $(,):\mathcal{H}^{\bullet}_R\times_M \mathcal{H}^{\bullet}_R 
\to 
CH,$ which restricts to a Hilbert $CH$-product on each fiber.
Further, on these bundles, Hilbert $C^*$-atlases have to be fixed 
such that their transition functions are maps into
$D^{\bullet}=\mbox{Aut}_{CH}(H^{\bullet}),$ i.e., the set of all homeomorphisms 
commuting with the action of $CH$ on $H^{\bullet}.$
 
\bigskip

We define a Hilbert $C^*$-bundle atlas on $\mathcal{H}^{0}_R$ by requiring that it is 
 the  
maximal smooth atlas satisfying the following conditions
\begin{itemize}
\item[1)] it contains $(M,J)$ and 
\item[2)] its transition functions are maps into $D^0$. 
\end{itemize}
Note that $D^0=\mbox{Aut}_{CH}(H^0) \simeq \mathbb{C}^{\times}$, which is quite 
restrictive  but it gives a clear picture of the situation.
Namely, the atlas equals $$\mathcal{A}= \{(U,\mu J_{|(p_0)^{-1}(U)})| \,     
\mu \in \mathbb{C}^{\times},\, U \mbox{ open in } M \}$$
This atlas is a subset of the  Kuiper structure on $\mathcal{H}
^0$ induced by $J.$
On $\mathcal{H}^{k}$ we consider  the tensor product atlases $\mathcal{A}_k,$  
which are the tensor products of the $C^{\infty}$-structure on 
$\bigwedge^kT^*M$ and of $\mathcal{A}.$ We  call them the {\it 
Hilbert $C^*$-bundle structure induced by $J$}.

\bigskip

{\bf Definition:} For a symplectic manifold $(M,\omega)$ possessing a 
metaplectic structure $(\pi_P,\Lambda)$, a compatible positive almost complex 
structure, a Kuiper map $J: \mathcal{H} \to M \times H,$ $m\in M,$
$f, h \in (\mathcal{H}^{\bullet}_R)_m,$ and $a \in CH,$ we set
\begin{align*}
&h \cdot a =   h \cdot_{\mathcal{H}} J_{CH}^{-1}(m,a) \hspace{.5cm} \mbox{and} 
\\
&(f,h) =\mbox{pr}_2^{CH} \left( J_{CH}(f,h)_{\mathcal{CH}}\right)
\end{align*}
where $\mbox{pr}_2^{CH}: M 
\times CH \to 
CH$ denotes the projection onto the second factor.

\bigskip

{\bf Remark:} The symbols $\cdot$ and $(,)$ should not be confused with the 
ones 
referring to  the action $H^{\bullet} \times CH \to H^{\bullet}$ and the 
$CH$-product $H^{\bullet} \times H^{\bullet} \to CH.$

\bigskip

Let us introduce the following operations. For $r \in \mathcal{P},$ $a,b \in 
CH,$ $v \in \check{H}$ 
 and  $f \in H,$ we set
\begin{align*}
 {\circ}_{\mathcal{CH}}: \mathcal{CH} \times_M \mathcal{CH} \to \mathcal{CH}& 
\quad [(r,a)] {\circ}_{\mathcal{CH}} [(r,b)] =   [(r,a  b)]\\
{\circ}_{CH}: (M \times CH) \times_M (M \times CH) \to M \times CH  & \quad 
(m,a)  {\circ}_{CH} (m,b) = (m,ab) \\
\otimes_{CH}: (M \times \check{H}) \times_M (M \times H) \to  M \times 
CH &\quad
(m,v) \otimes_{CH} (m,f)  = (m, v \otimes f)
\end{align*}
These maps  are well defined by Lemma 3.

\bigskip

{\bf Lemma 9:} Let $a,b \in CH,$  $f\in H$ and $m \in M.$ Then
\begin{align*}
&J^{-1}_{CH}(m,a\circ b)=J^{-1}_{CH}(m,a) \circ_{\mathcal{CH}} J_{CH}^{-1}(m,b) 
\mbox{ 
and }\\
&J^{-1}(m,f) \cdot_{\mathcal{H}} J^{-1}_{CH}(m,a) = J^{-1}(m,f \circ a)
\end{align*}

{\it Proof.} For any $q\in \mathcal{P},$ $v_a, v_b \in \check{H},$ and $f_a, 
f_b 
\in 
H,$ 
$([(q,v_a)]\otimes_{\mathcal{CH}}[(q,f_a)])\circ_{\mathcal{CH}}([(q,v_b)]
\otimes_{
\mathcal { CH } } [ (q , f_b) ]) = [(q, (v_a\otimes f_a) \circ (v_b\otimes 
f_b))] = [(q,f_a(v_b)v_a \otimes f_b)].$
Thus for $v,w \in \check{\mathcal{H}}_m$ and $f,h \in \mathcal{H}_m$
\begin{eqnarray}
      (v \otimes_{\mathcal{CH}} f) \circ_{\mathcal{CH}} 
(w\otimes_{\mathcal{CH}} 
h) 
= f(w)(v \otimes_{\mathcal{CH}} h) \label{blbost}
      \end{eqnarray}

 For $a=v_a \otimes f_a \in CH$ and  $b= v_b \otimes f_b \in CH$
\begin{align*} 
J_{CH}^{-1}(m,a \circ b) &= J_{CH}^{-1}\left(m, (v_a\otimes f_a) \circ (v_b 
\otimes f_b)\right)= 
f_a(v_b) J_{CH}^{-1}(m, v_a \otimes f_b)
\end{align*}
Since $J$ is a fiber-wise isometry
\begin{align*}
f_a(v_b) J_{CH}^{-1}(m, v_a &\otimes f_b)=f_a(v_b) \check{J}^{-1}(m,v_a) 
\otimes_{\mathcal{CH}}  J^{-1}(m,f_b)\\
&= (J^{-1}(m,f_a),J^{-1}((m,v_b)^{\flat_M}))_S  
\check{J}^{-1}(m,v_a)\otimes_{\mathcal{CH}} J^{-1}(m,f_b)\\
&=(J^{-1}(m,f_a)[\check{J}^{-1} 
(m,v_b)])\check{J}^{-1}(m,v_a)\otimes_{\mathcal{CH}}J^{-1}(m,f_b)
\end{align*}
Using (\ref{blbost}), we have  
$\left(J^{-1}(m,f_a)[\check{J}^{-1}(m,v_b)]\right)\check{J}^{-1}(m,v_a) 
\otimes_{\mathcal{CH}} J(m,f_b)$
\begin{align*}
&=(\check{J}^{-1}(m,v_a) \otimes_{\mathcal{CH}} J^{-1}(m,f_a)) 
\circ_{\mathcal{CH}} (\check{J}^{-1}(m,v_b) \otimes_{\mathcal{CH}} 
J^{-1}(m,f_b))\\
&=J_{CH}^{-1}(m, v_a \otimes f_a) \circ_{\mathcal{CH}} J_{CH}^{-1}(m, v_b 
\otimes 
f_b)=J_{CH}^{-1}(m,a) \circ_{\mathcal{CH}} J_{CH}^{-1}(m,b)
\end{align*}
 Summing-up, $J_{CH}^{-1}(m,a\circ b)=J_{CH}^{-1}(m,a) \circ_{\mathcal{CH}} 
J_{CH}^{-1}(m,b)$ 
for rank 1 operators $a$ and $b.$
For a linear combination of  finite rank maps, the result follows by linearity.
For a product of compact operators, one proceeds by taking  limits which is 
justified by noting that $J_{CH}$ is a homeomorphism (Sect. 3.3, proof of Thm. 
7).

 The second equation is proved in a similar way.

\hfill\(\Box\)

\bigskip 

{\bf Theorem 10:} Let $(M,\omega)$ be a symplectic manifold admitting a 
metaplectic structure $(\pi_P, \Lambda)$ and $J$ be a Kuiper map.
Then the reduced higher oscillator bundle  $\mathcal{H}^{\bullet}_R \to M$ with 
the Hilbert $C^*$-bundle structure  induced by a Kuiper map
$J$ is a smooth Hilbert $CH$-bundle for any compatible positive almost complex 
structure. 

{\it Proof.}  
First, we have to verify that the map $\cdot: \mathcal{H} \times CH \to 
\mathcal{H}$ is 
an action. Recall that according to the definitions of the 
$\mathcal{CH}$-action and the composition $\circ_{\mathcal{CH}}$, we have
$h \cdot_{\mathcal{H}} (a \circ_{\mathcal{CH}} b) = (h \cdot_{\mathcal{H}} a)
\cdot_{\mathcal{H}} b,$ whenever these operations make sense. 
For $m \in M,$ $h\in (\mathcal{H}^{\bullet}_R)_m$ and $a,b \in CH,$ we get 
according to Lemma 9
\begin{align*}
 h \cdot (ab) &= h \cdot_{\mathcal{H}} J^{-1}_{CH}(m,ab) = h 
\cdot_{\mathcal{H}} 
\left(J^{-1}_{CH}(m,a) \circ_{\mathcal{CH}} J^{-1}_{CH}(m,b)\right) \\
&= \left(h \cdot_{\mathcal{H}} J^{-1}_{CH}(m,a)\right) \cdot_{\mathcal{H}} 
J^{-1}_{CH}(m,b) = \left(h 
\cdot a\right) \cdot b
\end{align*}

Second, using Lemma 3, it is easy to verify that
$$(f,h \cdot_{\mathcal{H}} b)_{\mathcal{CH}} = (f,h)_{\mathcal{CH}} 
\circ_{\mathcal{CH}} b$$
for suitable $h,f \in \mathcal{H}_R^{\bullet}$ and $b \in \mathcal{CH}_R.$ For 
such elements and $a\in CH$ we get, using Lemma 9
\begin{align*}
(f, h \cdot a) &= \mbox{pr}_2^{CH}\left( J_{CH}\left( f, h \cdot_{\mathcal{H}} 
J_{CH}^{-1}(m,a)\right) \right)\\
\quad &=\mbox{pr}_2^{CH}\left( J_{CH}\left( 
(f,h)_{\mathcal{CH}} 
\circ_{\mathcal{CH}} J_{CH}^{-1}(m,a)\right)  \right)\\
\quad 
&=\mbox{pr}_2^{CH}\left(J_{CH}\left[J_{CH}^{-1}\left(J_{CH}(f,h)_{\mathcal{CH}}
\right) 
\circ_{\mathcal{CH}} J_{CH}^{-1}(m,a)\right]\right)\\
\quad &=\mbox{pr}_2^{CH}\left(J_{CH} \left[J_{CH}^{-1}\left(m, 
\mbox{pr}_2^{CH}\left(J_{CH}\left(f,h\right)_\mathcal{CH}\right)\circ 
a\right)\right] 
\right)\\
\quad &= \mbox{pr}_2^{CH}\left(m,\mbox{pr}_2^{CH}(J_{CH}(f,h)_{\mathcal{CH}})\circ 
a\right)\\
\quad &=\mbox{pr}_2^{CH}\left(J_{CH}\left(f,h\right)_{\mathcal{CH}}\right) 
\circ 
a = 
(f,h) \circ a
\end{align*}

The Hilbert $C^*$-bundle structure induced by $J$ on $\mathcal{H}^k,$ 
$k=0,\ldots, 
2n,$ is a  maximal set of smooth Hilbert $C^*$-bundle atlases smoothly compatible with 
$(M,J)$. Let us verify 
that the $CH$-action 
is 
smooth. The $CH$-action $\mathcal{H} \times CH \to \mathcal{H}$ in the map 
$(M,J)$ is given by
$(m,f) \mapsto (m, f\circ a),$ where $m\in M,$ $f\in 
(H_R^{\bullet})_m$ and $a\in CH.$  This map is infinitely many times 
Fr\'echet differentiable. Similarly, the $CH$-product is given by 
$((m,f),(m,h)) \mapsto (f,h) \in CH$ in $(M,J),$ $m\in M$ 
and  $f,h \in (H_R^{\bullet})_m,$ which is thus Fr\'echet 
smooth. 
\hfill\(\Box\)

\bigskip

{\bf Hilbert $C^*$-bundle structure on $\underline{H}$:} Let $g$ be the 
(Hodge type) Riemannian metric on the fibers of $\bigwedge^k T^*M$ induced 
by $\omega$ and a chosen compatible positive almost complex structure.
The tensor product bundle $\bigwedge^{\bullet}T^*M \otimes \underline{H} \to M$ 
of $\bigwedge^{\bullet}T^*M \to M$ and of the trivial
bundle $\underline{H} = M \times H \to M$ is a Hilbert 
$C^*$-bundle by setting $(\alpha \otimes (m,f))\cdot_M a = \alpha \otimes 
(m,f\cdot a) \in \bigwedge^{\bullet}T^*M \otimes 
\underline{H}$  and $(\alpha \otimes (m,f), \beta \otimes (m,h))_M = 
g(\alpha,\beta)(f,h)  \in CH,$ where 
$m\in M,$ $\alpha,\beta \in \bigwedge^{\bullet} T^*_mM,$ $f,h \in 
H$ and $a\in CH.$ For the smooth Hilbert $C^*$-bundle structure on 
$\bigwedge^{\bullet}T^*M \otimes \underline{H}$, we take the product of  
the maximal atlas of $\bigwedge^{\bullet}T^*M$ with the atlas   
$\{\left(U,\mu \mbox{Id}_{U \times H}\right)| \, \mu 
\in \mathbb{C}^{\times},\, U \mbox{ open in } M\}$ of $\underline{H}.$

\bigskip

{\bf Lemma 11:} Let $(M,\omega)$ be a symplectic manifold admitting a 
metaplectic structure $(\pi_P, \Lambda).$ For any
compatible positive almost complex structure and a Kuiper structure $J,$ the 
Kuiper map extends to a unitary isomorphism 
$\widetilde{J}: \mathcal{H}^{\bullet}_R \to \bigwedge^{\bullet}T^*M \otimes 
\underline{H}$ of Hilbert 
$CH$-bundles.

{\it Proof.} Obviously, $\mathcal{H}_R^k \simeq \bigwedge^kT^*M \otimes 
\mathcal{H}_R^0.$
We set $\widetilde{J}(\alpha \otimes h)= \alpha \otimes J(h),$ $\alpha  \in 
\bigwedge^kT^*M$ and $h\in \mathcal{H}^0_R,$ and extend it linearly.
Since $CH$ does not act on the form part, it is sufficient to verify that $J(h 
\cdot a) =J(h) \cdot_M a$ for any $h \in \mathcal{H}^0_R$ and $a\in CH.$
Representing $h=J^{-1}(m,f)$ and using Lemma 9
\begin{align*}
J(h\cdot a) &= J(J^{-1}(m,f) \cdot_{\mathcal{H}} J^{-1}(m,a)) = J(J^{-1}(m,h 
\circ a))\\
&= (m,h \circ a) = (m,h) \cdot_M a.
\end{align*}

It is easy to see
\begin{eqnarray}
(f,h)_{\mathcal{CH}} = f^{\sharp} \otimes_{\mathcal{CH}} h     \label{vzor}
\end{eqnarray}
 for $f, h \in \mathcal{H}^0_R.$
Indeed, representing $f, h$ by equivalence classes, we obtain \\
 $([(q,v)],[(q,w)])_{\mathcal{CH}} = [(q, v^{\sharp} \otimes w)] = f^{\sharp} 
\otimes_{\mathcal{CH}} h.$
  
 When  proving the unitarity of $\widetilde{J}$, we choose $J$-preimages of 
$f,h 
\in \mathcal{H}^0_R$ and use (\ref{vzor}), proceeding in a similar way 
 as in the case of $CH$-equivariance. The unitarity of $\widetilde{J}$ follows 
then from the fact that $CH$ does not act on the  wedge form part of
$\bigwedge^{\bullet}T^*M \otimes \underline{H}.$  
\hfill\(\Box\)

\subsection{Associated elliptic complex to a Kuiper connection}

We introduce the complex and prove that it is a complex of 
$C^*$-operators. Further we investigate its cohomology groups. For that we 
focus to the topology of section spaces of the higher oscillator bundles 
and the topology of their tensor products and quotients.

Let $(M,\omega)$ be a symplectic manifold admitting a metaplectic structure 
$(\pi_P, \Lambda)$ and $J: \mathcal{H} \to M \times H$ be a 
normalized Kuiper map. 
The $T(M\times H)$-valued  $1$-form $\Phi^H$ defining the canonical 
product  connection on $M \times H 
\to M$ is 
given by $$\Phi^H(t_m, v_h)=(0_m,v_h)$$
where  $m\in M,$ $t_m \in  T_mM,$ $h\in H,$  $v_h \in T_hH,$ and $0_m$ is the 
zero tangent vector in $T_mM.$ Note that we identify $T(M\times H) \simeq TM 
\times 
TH$.
The Kuiper map gives rise to a connection $\Phi$ in $\mathcal{H}$ given by 
the prescription $\Phi=TJ^{-1} \circ \Phi^H \circ TJ$ which we call the 
{\it Kuiper connection induced by} $J,$ expressing its roots.    

\bigskip

 Let us notice that a 
different  description of a connection on symplectic spinors is chosen by 
Herczeg, Waldron \cite{HerczegWaldron}  in the 
setting of quantum physics for contact manifolds. This connection is used
already in Kr\'ysl \cite{KryslDGA0}.

\bigskip

{\bf Remark:} 
The distribution  $$\mathcal{H} \ni f=[(r,h)] \mapsto \mathcal{L}_f = 
\left(T_{(\pi_P(r),h)}J^{-1}\right)(T_mM,0)\subseteq 
T_f\mathcal{H}$$ where  
$r \in \mathcal{P}, h \in H,  m=\pi_P(r)$ and $0 \in  T_hH$ (the zero vector 
of $T_hH$) is the horizontal  distribution defined by $\Phi.$ In this way, we 
get 
a splitting of the vector bundle $T\mathcal{H} = 
\mathcal{L} \oplus \textrm{Ker}\,(Tp_0).$  The vertical part of the bundle 
$\mathcal{H}$ is denoted by $\mathcal{V}.$

\bigskip

{\bf Gateaux and Fr\'echet differentials in  geometric setting:} We  
introduce a notation and recall some facts from analysis on Banach spaces.
Let $(Df)(h,v)$ denote the  Gateaux derivative of $f$ at $h$ in direction $v.$ 
For each $h \in H,$ $T_hH = \{v_h, v \in H\},$ where $v_hf=(Df)(h,v)$ for any
smooth function $f$ defined in a neighborhood of $h \in H.$   
For any point $h \in H,$ let us define $i_h:T_hH \to H$ by $i_h(v_h)=v,$ $v \in 
H.$ It is easy to see that 
$i_h$ is an isomorphism of vector spaces.

\bigskip

{\bf Vertical lifts, connectors and induced connections:} We recall 
basic principles for handling of associated 
connections. It is  based on the  approach described in Michor 
\cite{Michor} and Kol\'a\v{r}, Michor, Slov\'ak \cite{KMS}.
For a vector field $X\in \mathfrak{X}(M)(=\Gamma(TM)),$ the connection on 
$\mathcal{H}$ induces a 
covariant derivative $\nabla_X: \Gamma(\mathcal{H}) \to \Gamma(\mathcal{H})$
defined by $\nabla_X s = \kappa \circ Ts \circ X,$ where
$X: M \to TM$ is a vector field on $M,$ $Ts: TM \to T\mathcal{H}$ is the tangent
map of a bundle section $s: M \to \mathcal{H},$ and
$\kappa: T\mathcal{H} \to \mathcal{H}$ is  the connector of $\Phi.$
We have also a covariant derivative $\nabla^H_X$ acting on sections of $M\times 
H \to M$ which is given by $\nabla^H_X s = \kappa^H \circ Ts \circ X,$ where 
$s: M \to M \times H$ and $\kappa^H$ is the connector of $\Phi^H.$
Let us denote the vertical lifts of $(H\times M \to M,\Phi^H)$ and 
$(\mathcal{H} 
\to M, \Phi)$ by $vl^H$
and $vl,$ respectively. Recall that $vl: \mathcal{H} \times_M \mathcal{H} \to 
\mathcal{V}$ and similarly 
$vl^H: (M \times H) \times_M (M \times H) \to TM \times H.$
We denote the projections $\mbox{pr}_2: \mathcal{H} \times_M 
\mathcal{H} \to \mathcal{H}$ and $\mbox{pr}_2^H:(M\times H)\times_M (M\times H) 
\to M\times H$
onto the respective second factors. With this notation, the connector $\kappa$ 
is the map $\mbox{pr}_2 \circ vl^{-1} \circ \Phi$ and similarly for $\kappa^H.$ 

\bigskip

{\bf Lemma 12:} In the situation described in the previous  paragraph, we have 
for any $m\in M,$ and $h,v \in H$
\begin{itemize} 
\item[1)] $(vl^H)^{-1}(0_m,v_h)=\left((m,h),(m,v)\right) \in (M\times H) 
\times_M 
(M\times H)$
\item[2)] $vl = TJ^{-1} \circ vl^H \circ (J \times_M J)$
\item[3)] $\mbox{pr}_2^H \circ  (J \times_M J) = J \circ \mbox{pr}_2 $ 
\end{itemize}

{\it Proof.} For  $m \in M$ and
$h, v \in (M \times H)_m,$ 
$vl^H\left((m,h),(m,v)\right)=$\\
$=\frac{d}{dt}_{|t=0}\left( (m,h) 
 +_M t\odot_M (m,v)\right),$ where
$+_M$ and $\odot_M$ denote the addition and the  multiplication by 
scalars in 
$M \times H \to M,$ respectively, as introduced  in 3.3.
We have $\frac{d}{dt}_{|t=0} ((m,h)+_Mt\odot_M(m,v))=\frac{d}{dt}_{|t=0}(m, 
h+tv) = 
(0_m,v_{h}).$ 

For $x,y \in \mathcal{H}_m,$ $vl(x,y) = \frac{d}{dt}_{|t=0}(x+ty).$ Because 
$J$ is fiber-wise linear,   $vl(x,y)= 
 \frac{d}{dt}_{|t=0}J^{-1}J(x+ty)=
TJ^{-1}\frac{d}{dt}_{|t=0}(h +_M t\odot_M v) = TJ^{-1}vl^H(h,v),$ where 
$h = J(x)$ and $v=J(y).$

The third equality is obvious.
\hfill\(\Box\)

\bigskip

{\bf Lemma 13:} For $X \in \mathfrak{X}(M)$ and $s\in \Gamma(\mathcal{H})$ 
$$J \circ \nabla_X s = \nabla_X^H (J \circ s)$$

{\it Proof.} We use Lemma 12  to prove the identity. For $X\in \mathfrak{X}(M)$ 
and $s\in \Gamma(\mathcal{H})$
\begin{eqnarray*}
\nabla_X^H(J \circ s) &=& \kappa^H \circ TJ \circ Ts \circ X\\
&=& \mbox{pr}_2^H \circ {vl^H}^{-1} \circ \Phi^H \circ TJ \circ Ts \circ X\\
&=& \mbox{pr}_2^H \circ (J \times_M J) \circ vl^{-1}  \circ TJ^{-1} \circ 
\Phi^H 
\circ 
TJ \circ Ts \circ X\\
&=& \mbox{pr}_2^H \circ (J \times_M J) \circ vl^{-1} \circ \Phi \circ Ts \circ 
X\\
&=& J \circ \mbox{pr}_2 \circ vl^{-1}   \circ \Phi \circ Ts \circ X = J \circ 
\nabla_Xs
\end{eqnarray*}
\hfill\(\Box\)

\bigskip

{\bf Remark:} Since we shall compute derivations of the $CH$-action on 
$\underline{H}$ (by a fixed element), 
we  use a notation  which is more usual for maps.  
For any $a \in CH$ and $(m,v) \in M \times H,$ we set $G_a(m,v)=(m, v \circ a)$ 
and $g_a(v) = v \circ a,$
defining maps $G_a: M \times H \to M \times H$ and $g_a: H \to H.$ Note that 
$G_a = 
\mbox{Id}_M \times g_a.$
 
\bigskip

{\bf Lemma 14:} For any $X \in \mathfrak{X}(M),$ $s\in \Gamma(\underline{H})$ 
and 
$a \in CH$ 
$$\nabla^H_X(G_a \circ s) = G_a \circ \nabla^H_X s$$

{\it Proof.} 
For  $v,h\in H$ and $a \in CH,$ $(Dg_a)(h,v)=\lim_{t\to 0} [t^{-1}(h \circ a  
+ t v \circ a - h \circ a)] =  v \circ a$.
For $f\in C^{\infty}(H),$ $(Tg_a)(v_h)f = v_h(f \circ g_a) =  
D(f \circ g_a)(h,v)= $ $(Df)(h \circ a,$ $(Dg_a)(h,v)) = (Df)(h \circ a, v 
\circ 
a)$ by the chain rule. 
The last expression equals to $(v \circ a)_{h \circ a} f.$
Thus  $(Tg_a)v_h = (v \circ a)_{h \circ a},$ which implies
\begin{align}
TG_a(t_m,v_h) = \left(t_m, (Tg_a)(v_h)\right) = (t_m, (v \circ a)_{h \circ a}) 
\label{pomoc}
\end{align}
where $m \in M$ and $t_m \in T_mM.$

By Lemma 12 
\begin{eqnarray*}
(G_a \circ \kappa^H)(t_m,v_h)&=&(G_a \circ \mbox{pr}_2^H \circ{vl^H}^{-1} 
\circ \Phi^H)(t_m,v_h) \\
&=& (G_a \circ \mbox{pr}_2^H \circ {vl^H}^{-1})(0_m,v_h)\\
&=&(G_a \circ \mbox{pr}_2^H)((m,h),(m,v))=G_a(m,v) = (m,v \circ a)
\end{eqnarray*}
 
Further, we get by Lemma 12 and (\ref{pomoc})  that
\begin{eqnarray*}
(\kappa^H \circ TG_a)(t_m,v_h)&=& (\mbox{pr}_2^H \circ{vl^H}^{-1}\circ \Phi^H 
\circ TG_a)(t_m,v_h)\\
&=&(\mbox{pr}_2^H \circ{vl^H}^{-1} \circ\Phi^H)(t_m,i^{-1}_{h \circ a}(v \circ 
a))\\ 
&=& (\mbox{pr}_2^H \circ{vl^H}^{-1})(0_m, i^{-1}_{h \circ a}(v \circ a)) = (m,v 
\circ a)
\end{eqnarray*}

Summing-up, $G_a \circ \kappa^H = \kappa^H \circ TG_a.$
The theorem follows by 
$\nabla_X(G_a \circ s)= \kappa^H \circ TG_a \circ Ts \circ X = G_a \circ 
\kappa^H \circ Ts \circ X = G_a \circ \nabla^H_X s$ for any
  $s\in \Gamma(\underline{H})$ and $X \in \mathfrak{X}(M).$
\hfill\(\Box\)

\bigskip

{\bf Hilbert $C^*$-module structure on section spaces:}
Notice that $CH$ acts not only on the higher oscillator bundle but 
also on  its smooth 
sections space $\Gamma(\mathcal{H}^{\bullet}).$
Namely, $(s \cdot a)_m = s_m\cdot a$
for any $a\in CH,$  $s \in \Gamma(\mathcal{H}^{\bullet})$ and $m \in 
M.$\footnote{At the left-hand side a new action is defined.} In this 
way, 
$\Gamma(\mathcal{H}^{\bullet})$ is a right $CH$-module. The space of smooth 
sections of 
the reduced higher oscillator bundle $\mathcal{H}_R^{\bullet}$  is also a 
pre-Hilbert 
$CH$-module if $(M, \omega)$ is compact.
The $CH$-valued product is given by $$(s,s') = \int_{m \in M} (s_m,s'_m) 
\mbox{vol}_m$$
where $s,s' \in \Gamma(\mathcal{H}^{\bullet})$ and $\mbox{vol}_m$ denotes the 
volume form in $m$ 
induced by $\omega_m.$ The integral is the Pettis integral (for the measure on 
$M$ 
induced by the volume form $\mbox{vol}_m$) of
Pettis $CH$-valued integrable function on $M.$ (See Ryan 
\cite{Ryan} or Grothendieck \cite{Groth}.) 
Similarly, we consider the 
space  
$\Gamma(\underline{H})$ of smooth sections of the product bundle $\underline{H} 
\to M$ as a 
pre-Hilbert $CH$-module as well. Completions of this pre-Hilbert module
(see Fomenko, Mishchenko \cite{FM}) with respect to the induced norm make 
it a Hilbert $C^*$-module. (Let us notice that for reasons of analysis, several
further $C^*$-valued products of `Sobolev type' are introduced in \cite{FM} 
besides 
the $CH$-product given above. However, the integral seems to be rather not 
specified there.)

\bigskip

 {\bf Corollary 15:} For a  symplectic manifold $(M,\omega)$ admitting a 
metaplectic structure
$(\mathcal{P}, \pi_P),$  a Kuiper map $J$ and a vector field $X$ on $M,$ 
the operator 
$$\nabla_X:\Gamma(\mathcal{H}) \to \Gamma(\mathcal{H})$$ is a first order 
$CH$-equivariant differential operator.

{\it Proof.}
 The map $\nabla_X$ is a first order differential operator for any $X\in 
\mathfrak{X}(M)$ (all   operations involved in its definition are of zero
order, except of $s \mapsto Ts$).
Using Lemmas 9, 13 and 14  for $a \in CH,$ $X \in \mathfrak{X}(M)$ and 
$s \in \Gamma(\mathcal{H})$
\begin{eqnarray*}
\nabla_X (s \cdot a) &=& J^{-1}\circ  \nabla_X^H (J \circ (s \cdot a))= 
J^{-1}\circ  \nabla^H_X(J \circ J^{-1} \circ G_a \circ J \circ s)\\
&=& J^{-1}\circ  \nabla_X^H (G_a \circ  J \circ s) =(J^{-1} \circ G_a) \circ 
\nabla^H_X (J \circ s)\\
&=& (J^{-1} \circ G_a \circ J) \circ \nabla_X s  = (\nabla_X s) \cdot a
\end{eqnarray*}
\hfill\(\Box\)

\bigskip

{\bf Remark} (Connection to non-commutative geometry): It is  easy to see that 
$\nabla_X$ is $CH$-hermitian
in the sense that $X(s,s') = (\nabla_X s, s')+ (s,\nabla_X s')$ where $s,s' \in 
\Gamma(\mathcal{H}_R^{\bullet}).$ 
In Dubois-Violette, Michor  \cite{DM}, central $A$-bimodules and connections 
for them are investigated. To some extent, our structures fit into the specific 
{\it non-commutative geometry concept} described there.
Namely, for the not necessarily normed algebra $A$, considered in \cite{DM}, we 
take the smooth 
sections of the bundle of measurements,  $A=\Gamma(\mathcal{CH}),$ and for the 
module, we take the  smooth sections  $\Gamma(\mathcal{H}^{\bullet}).$ The 
$A$-valued metric 
$h$ (p. 229 in \cite{DM}) is then the `sectionised'  $\mathcal{CH}$-product,   
i.e.,  $h: \Gamma(\mathcal{H}^{\bullet}_R) \times 
\Gamma(\mathcal{H}_R^{\bullet}) \to 
\Gamma(\mathcal{CH}),$ $h(s,s')(m)=(s_m,s'_m)_{\mathcal{CH}}$ where 
$s,s' \in \Gamma(\mathcal{H}^{\bullet})$ and $m\in M.$ 
However, the 
vector fields $X$ in our paper are generically not all derivations of 
$A=\Gamma(\mathcal{CH}).$ (They coincide iff the manifold is discrete as a 
topological space.) The module $\Gamma(\mathcal{H}_R^{\bullet})$ with the 
$\Gamma(\mathcal{CH})$-hermitian connection  may be considered as a modification of
a non-commutative hermitian manifold as considered in the program of Connes \cite{Connes}. 
They have a global realization by topological and differential structures.

\bigskip

Let $(e_i)_{i=1}^{2n}$ and $(\epsilon^i)_{i=1}^{2n}$ be  dual local frames of 
$TM \to M$ and $T^*M \to M$ respectively. 
Recall that the  exterior covariant derivatives 
$d_k^{\Psi}:\Gamma(\bigwedge^kT^*M \otimes \mathcal{F}) \to 
\Gamma(\bigwedge^{k+1}T^*M \otimes \mathcal{F})$ induced by a vector bundle
connection $\Psi,$ defined on a vector bundle $\mathcal{F} \to M,$ are given by
   $$d^{\Psi}_k:  \alpha \otimes s \mapsto 
d_k\alpha \otimes s + \sum_{i=1}^{2n}\epsilon^i \wedge \alpha \otimes 
\nabla_{e_i}^{\Psi}s$$ 
where $\nabla^{\Psi}$ denotes the covariant derivative associated to $\Psi$ and 
$d_k$ denotes the de Rham differential.   They are  extended  linearly to 
non-homogeneous elements. 

\bigskip

{\bf Theorem 16:} The   exterior covariant  derivatives induced by the 
Kuiper connection $\Phi$ form 
a cochain 
complex in the category of right-modules over the ring $CH$ of compact 
operators, i.e.,
\begin{align*}
&d^{\Phi}_k(s \cdot a) = (d^{\Phi}_k s)\cdot a \, \mbox{  and}\\
&d_{k+1}^{\Phi} \circ d_k^{\Phi} = 0
\end{align*}
for any $a\in CH$ and $s\in \Gamma(\mathcal{H}^{\bullet}_R).$

{\it Proof.} 
The second formula follows  by  applying the definition of the exterior 
covariant derivatives and  using the fact that $\Phi$ is flat. 
 
Using  Corollary 15, for any $\phi = \alpha \otimes s \in 
\Gamma(\mathcal{H}^k)$ and $a\in CH$ 
\begin{eqnarray*}
d^{\Phi}_k \left((\alpha \otimes s) \cdot a\right) &=& d^{\Phi}_k\left(\alpha 
\otimes (s 
\cdot a)\right)\\
&=& d_k\alpha \otimes (s \cdot a) +  \sum_{i=1}^{2n}
\epsilon^i \wedge 
\alpha \otimes \nabla_{e^i}^{\Phi} (s\cdot a)\\
&=& (d_k\alpha \otimes s)\cdot a + \sum_{i=1}^{2n}\epsilon^i \wedge \alpha 
\otimes 
(\nabla_{e_i}^{\Phi} s) \cdot a\\
&=& \left( d_k\alpha \otimes s  + \sum_{i=1}^{2n} \epsilon^i \wedge \alpha 
\otimes \nabla_{e_i}^{\Phi} s \right) \cdot a\\
&=& \left(d^{\Phi}_k(\alpha \otimes s)\right)\cdot a
\end{eqnarray*}
\hfill\(\Box\)

\bigskip

{\bf Definition:} For a symplectic manifold admitting a metaplectic structure 
$(\pi_{P}, \Lambda),$  a compatible positive almost complex structure and a 
normalized Kuiper map $J,$ we  call the cochain complex $d^{\Phi}_{\bullet}= 
(\Gamma(\mathcal{H}_R^k), d^{\Phi}_k)_{k\in \mathbb{Z}}$ in the category of 
pre-Hilbert $CH$-modules
the {\it Kuiper complex induced by $J$}.

\bigskip 

{\bf Remark:} Let $A$ be a $C^*$-algebra. In  categories
of pre-Hilbert and Hilbert $A$-modules and adjointable maps between them, the 
images of morphisms need not have closed range. (See Solovyev, Troitsky 
\cite{ST}.)
Consequently, the quotient topology on the cohomology groups of complexes in  
 these categories may be non-Hausdorff and thus, the quotients may not be 
Hilbert $A$-modules in the quotient topology. (See  \cite{KryslJGP2} for simple 
examples.)

\bigskip

{\bf Convention on the cohomology groups and tensor products:} For a compact 
manifold $M,$ we denote the space of $k$-cocycles and the $k$-th  cohomology 
group of a cochain complex $D_{\bullet}$ on $\mathcal{H}^{\bullet}_R$ by 
$Z^k(D_{\bullet})$ and
$H^k(D_{\bullet}),$ respectively.  The topology on 
$\Gamma(\mathcal{H}^{\bullet}_R)$
is determined   by the 
norm induced from the pre-Hilbert $C^*$-module structure of 
$\Gamma(\mathcal{H}^{\bullet}_R).$ For the topology on the cohomology groups, 
we take 
the quotient topology. In the case of the de Rham complex on $M$, we denote the 
space of $k$-cocycles by $Z^k_{dR}(M),$ the space of $k$-coboundaries by
$B^k_{dR}(M)$ and the cohomology group by 
$H^k_{dR}(M),$ and 
consider the  quotient topology on it as well. On the space of exterior 
differential $k$-forms $\Omega^k(M)$ on the manifold $M$, we take the norm 
topology given by the  Hodge scalar  product.
When we form a tensor product of pre-Hilbert spaces, we mean 
the algebraic tensor product equipped with the scalar product induced by the 
scalar products on the individual factors (especially, no 
completions  involved).

\bigskip

{\bf Lemma 17:} Let $(h_i)_{i\in \mathbb{N}}$ be a complete orthonormal system 
on $H$ and, for $i \in \mathbb{N},$ let $\underline{h}_i$ be the constant 
extension of $h_i$ to a section of $\Gamma(\underline{H}).$  If  $h= 
\sum_{i=1}^{\infty} c_i\underline{h}_i$ is a $C^1$-section of $\underline{H}$, 
then 
for any $i \in \mathbb{N},$ $c_i$ are $C^1$-functions and for any $v\in TM$ 
$$\nabla_v^H h = 
\sum_{i=1}^{\infty} v(c_i) \underline{h}_i$$

{\it Proof.} Since the pointwise scalar product $(,)$ and $h$ are 
continuous, we see 
that $c_j$ is continuous for any $j\in \mathbb{N}$ by taking the scalar product 
of $h$ and $\underline{h}_j.$ 
 Let us choose a map $(U,\phi),$ $x\in U,$ $v\in T_xM$ and compute `in coordinates'
\begin{eqnarray*}
\nabla_v h&=&\lim_{t\to 0} \frac{1}{t} \sum_{i=1}^{\infty} 
(c_i(x+tv)\underline{h}_i(x+tv) - c_i(x)\underline{h}_i(x)) \\
&=&\lim_{t \to 0} \sum_{i=1}^{\infty} \frac{1}{t}(c_i(x+tv)-c_i(x)) 
\underline{h}_i(x).  
\end{eqnarray*}
 Taking the inner product with $\underline{h}_j(x),$ we get 
\begin{eqnarray*}
(\nabla_v h, \underline{h}_j(x)) &=& (\lim_{t \to 0} \sum_{i=1}^{\infty} 
\frac{1}{t}(c_i(x+tv)-c_i(x)) \underline{h}_i(x), \underline{h}_j(x))
\end{eqnarray*}
Since the pointwise inner product $(,),$ $c_i$'s and $\underline{h}_i$'s are continuous, we can take the limit out of the bracket getting
\begin{eqnarray*}
&&\lim_{t \to 0} \frac{1}{t}(\sum_{i=1}^{\infty} (c_i(x+tv)-c_i(x))
\underline{h}_i(x), \underline{h}_j(x)) = \lim_{t\to 0} 
\frac{1}{t}(c_j(x+tv)-c_j(x))\\ 
&&= v(c_j)
\end{eqnarray*}
Especially, the derivation of $c_j$ in direction $v$ exists and it is 
continuous since the scalar product and $\underline{h}_j$ are continuous 
  and since $\nabla_v h$ is continuous by the assumption.
From the uniqueness of  coordinates in a Hilbert space, we get  that
$\nabla_v h = \sum_{i=1}^{\infty} v(c_i)\underline{h}_i.$    
\hfill\(\Box\)

\bigskip

{\bf Theorem 18:} Let $(M,\omega)$ be a compact symplectic manifold of 
dimension $2n$ 
admitting a metaplectic structure $(\pi_P,\Lambda)$ and $J$ be a 
normalized Kuiper map.
Then for any $k \in \mathbb{Z}$ and any compatible positive almost complex 
structure,  the cohomology groups 
$H^k(d_{\bullet}^{\Phi})$ of the Kuiper complex induced by $J$
have a structure of  finitely generated projective Hilbert $CH$-modules whose 
norm topology coincides with the quotient topology.

{\it Proof.} 
1) By Lemma 11, the bundles $\mathcal{H}^k_R$ and $\bigwedge^kT^*M \otimes 
\underline{H}$ are isomorphic
as smooth Hilbert $C^*$-bundles. 
  The  map $$j: \Gamma(\mathcal{H}^{\bullet}_R) \to  
\Gamma(\bigwedge^{\bullet}T^*M \otimes \underline{H}) \mbox{ given by } 
j(s) = \widetilde{J} \circ s$$  where $s \in 
\Gamma(\mathcal{H}_R^{\bullet}),$ 
preserves the $CH$-products and it is a $CH$-homomorphism as follows from  
Lemma 11 and the linearity of the Pettis integral with respect to continuous 
maps (see, e.g., \cite{Ryan}), respectively.
The same is true for the inverse 
$j^{-1}(s) = \widetilde{J}^{-1} \circ s$. In particular,
$\Gamma(\bigwedge^k T^*M \otimes \underline{H})$   and 
$\Gamma(\mathcal{H}_R^{\bullet})$ are isomorphic as pre-Hilbert $CH$-modules.

2) For $\alpha \otimes s \in \Gamma(\mathcal{H}_R^k),$ we have using Lemma 13
\begin{align*}
&d_k^{\Phi^H} \left( j (\alpha \otimes s) \right) =d_k^{\Phi^H} \left(\alpha 
\otimes (J \circ s)\right) = d_k\alpha \otimes J \circ s + \sum_{i=1}^{2n} 
\epsilon^i \wedge \alpha \otimes \nabla_{e_i}^H
(J \circ s) \\
&= d_k\alpha \otimes J \circ s + \sum_{i=1}^{2n} \epsilon^i \wedge \alpha \otimes J 
\circ \nabla_{e_i} s = \widetilde{J} \circ d_k^{\Phi}(\alpha \otimes 
s)=j\left(d_k^{\Phi}(\alpha \otimes s)\right)
\end{align*}
where $(e_i)_{i=1}^{2n}$ and $(\epsilon)_{i=1}^{2n}$ are dual local frames of 
$M.$ Thus $j$ is a cochain map.
Because $j$ is a pre-Hilbert $CH$-isomorphism (item 1), 
the corresponding cohomology 
groups $H^k(d_{\bullet}^{\Phi})$ and $H^k(d_{\bullet}^{\Phi^H})$ of 
$(\Gamma(\mathcal{H}_R^k), d_k^{\Phi})_k $ and $(\Gamma(\bigwedge^kT^*M 
\otimes \underline{H}), d_k^{\Phi^H})_k$   are isomorphic as 
 $CH$-modules and topological vector spaces when considered with  the
quotient topologies.

3) Note that $H^k_{dR}(M) \otimes H$ is a $CH$-module by defining 
$([\alpha]\otimes h)\cdot a = [\alpha] \otimes (h \cdot a)$ for $\alpha \in 
\Omega^k(M),$ $h\in  H$ and $a\in CH.$ We show that 
the cohomology groups of  
$(\Gamma(\bigwedge^kT^*M \otimes 
\underline{H}), d_{k}^{\Phi^H})_k$ are isomorphic to $H^k_{dR}(M) \otimes H$ as 
$CH$-modules. Let us consider the linear 
extension, denoted by $[\gamma],$ of the map 
$H^k_{dR}(M) \otimes H \ni [\alpha]\otimes h \mapsto [\alpha \otimes 
\underline{h}] 
\in H^k(d_{\bullet}^{\Phi^H})$
where $\underline{h}$ is the constant extension of an element $h \in H.$
It is easy to check that the resulting map is well defined (maps cocycles into 
cocycles and vanishes on coboundaries).
 Let us analyze the continuity of $[\gamma].$  For it, we denote   the 
extension map $h \mapsto \underline{h}$  by $\mbox{ext}.$ Then 
$[\gamma]$ fits 
into the following commutative diagram 
$$\begin{xy}\xymatrix{
Z^k_{dR}(M)\otimes H \ar[d]^{q_1 \otimes \mbox{Id}_H} \ar@{->}[r]^-{\mbox{Id}
 \otimes 
ext} 
&Z^k(d_{\bullet}^{\Phi^H})\ar[d]^{q_2}\\ 
H^k_{dR}(M)\otimes H \ar@{->}[r]^-{[\gamma]}
 &H^k(d_{\bullet}^{\Phi^H})}\end{xy}$$
where $q_i,$ $i=1,2,$ denote the appropriate quotient projections and 
$\mbox{Id}$ 
denotes $\mbox{Id}_{Z^k_{dR}(M)}$.
Since $q_2,$ $\mbox{Id}$ and $\mbox{ext}$ are continuous and 
$q_1\otimes \mbox{Id}_H$ is open, $[\gamma]$ is continuous. 

For to prove that  $[\gamma]$ is bijective, 
let us choose a complete orthonormal system $(h_i)_{i \in \mathbb{N}}$ of $H$ 
and let us consider an element $f=[\sum_{j=1}^r \alpha_j \otimes 
(\sum_{i=1}^{\infty}c_{ij}\underline{h_i})]$ of $H^k(d_{\bullet}^{\Phi^H}),$ 
where
$c_{ij}$ are functions on $M.$ The element $f'=\sum_{i=1}^{\infty}[\sum_{j=1}^r 
c_{ij}\alpha_j]\otimes h_i \in H^k_{dR}(M)\otimes H$ is a well defined 
element of the tensor product, as can 
be seen by the following consideration. Since  
$c_{ij}$ are smooth (by an iterative use of Lemma 17),
$\sum_{j=1}^rc_{ij}\alpha_j$ is smooth for each $i \in \mathbb{N}.$ Further, the 
 sum $\sum_{i=1}^{\infty}\sum_{j=1}^r 
c_{ij}\alpha_j\otimes h_i,$ denoted by $f'',$ converges pointwise on the
manifold  to
$\sum_{j=1}^r \alpha_j\otimes \sum_{i=1}^{\infty}c_{ij}\underline{h_i},$ that converges 
pointwise   in each point on the manifold by the 
choice of the representation for $f$.
  Let us consider the operator 
$T: Z^k_{dR}(M) \otimes H \to H^k_{dR}(M) \otimes H$ given as the 
linear extension of $T(\alpha \otimes 
h) = [\alpha] \otimes h.$ Let us see that this  operator is continuous. 
 Since $M$ is compact, $H^k_{dR}(M)$ is finite dimensional and $B^k_{dR}(M)$ is 
orthogonally complemented by the usual Hodge theory. Therefore, we can define the 
orthogonal projection $\mbox{pr}$ of $Z^k_{dR}(M)$ 
onto the complement, which is the space of harmonic forms.
The quotient norm of an element $[\alpha]$ in the de Rham cohomology group 
equals to the norm of $\mbox{pr}(\alpha)$
in the cocycle space. Consequently, the norm of $T$ is bounded by the norm of 
$\mbox{pr}$ which is one. Now, $f'$ is the $T$-image  of $f''$ by the 
continuity of $T.$ Especially, $f'$ is  well defined (from the point of view
of convergence).
The facts that  $\sum_{j=1}^r \alpha_j \otimes 
\sum_{i=1}^{\infty}c_{ij}\underline{h_i}$ is closed and that $\underline{h}_i$'s are 
constant imply 
\begin{eqnarray*}
0 &=& d_{k}^{\Phi^H}(\sum_{j=1}^r \alpha_j \otimes 
\sum_{i=1}^{\infty}c_{ij}\underline{h_i}) \\
&=&\sum_{j=1}^r \left(d_{k}\alpha_j \otimes 
\sum_{i=1}^{\infty} c_{ij}\underline{h}_i + \sum_{p=1}^{2n}\epsilon^p \wedge 
\alpha_j \otimes \sum_{i=1}^{\infty} e_p(c_{ij})\underline{h}_i\right)\\
&=&\sum_{j=1}^r\left(d_{k}\alpha_j \otimes\sum_{i=1}^{\infty} 
c_{ij}\underline{h}_i + \sum_{i=1}^{\infty}(d_0c_{ij})\wedge
\alpha_j \otimes \underline{h}_i\right)\\
&=& \sum_{i=1}^{\infty}\sum_{j=1}^rd_k(c_{ij}\alpha_j) \otimes \underline{h}_i
\end{eqnarray*}
Since $(h_i)_{i\in \mathbb{N}}$ is linearly independent,   
$\sum_{j=1}^rc_{ij}\alpha_j$ is closed for each $i\in \mathbb{N}.$ Especially,
$f'\in H^k_{dR}(M)\otimes H.$ 
Since $[\gamma]$ is continuous, 
$f'=\sum_{i=1}^{\infty}[\sum_{j=1}^r 
c_{ij}\alpha_j]\otimes h_i$ is a preimage of $f.$ Summing-up, $[\gamma]$ is surjective. 

For the injectivity, let us 
suppose that there is an element   $\sum_{l=1}^s\sum_{i=1}^{\infty} 
[\alpha_l] \otimes c_{li}h_i$ mapped to $[\sum_{l=1}^s 
\sum_{i=1}^{\infty}\alpha_l \otimes c_{li}\underline{h}_i]=0$ by $[\gamma],$ 
where $c_{li}$ are complex numbers.  
Thus, there exist smooth  exterior differential $(k-1)$-forms $\widetilde{\alpha}_j,$ 
smooth functions $\widetilde{c}_j$ and smooth sections 
$v_j \in \Gamma(\underline{H}),$  $j=1, \ldots, r,$ such that 
\begin{eqnarray*}
\sum_{l=1}^s\sum_{i=1}^{\infty}c_{li}\alpha_l \otimes \underline{h}_i = 
d^{\Phi^H}_{k-1} \left(\sum_{j=1}^r\widetilde{c}_j\widetilde{\alpha}_j \otimes 
v_j\right) 
\end{eqnarray*} For any $v_j,$ there are
smooth functions $\widetilde{c}_{ij}$ on $M$ (Lemma 17) for which 
$v_j=\sum_{i=1}^{\infty}\widetilde{c}_{ij}\underline{h}_i.$
Using the Leibniz rule, we get
\begin{eqnarray*}
&&\sum_{l=1}^s\sum_{i=1}^{\infty}c_{li} \alpha_l \otimes 
\underline{h}_i=d_{k-1}^{\Phi^H}\left(\sum_{j=1}^r \widetilde{c}_j\widetilde{\alpha}_j 
\otimes \sum_{i=1}^{\infty}\widetilde{c}_{ij}\underline{h}_i\right)=\\ 
&&d_{k-1} \left(\sum_{j=1}^r \widetilde{c}_j\widetilde{\alpha}_j\right) \otimes  \sum_{i=1}^{\infty}\widetilde{c}_{ij}\underline{h}_i
+ \sum_{p=1}^{2n}\sum_{j=1}^r \widetilde{c}_j\epsilon^p \wedge 
\widetilde{\alpha}_j \otimes \nabla^H_{e_p}(\sum_{i=1}^{\infty} 
\widetilde{c}_{ij} \underline{h}_i) =:(*)
\end{eqnarray*}
 
Using Lemma 17, we can take the derivative behind the infinite sum, getting 
\begin{eqnarray*}
(*) &=& \sum_{j=1}^r \left( d_{k-1}(c_j  \widetilde{\alpha}_j) \otimes  \sum_{i=1}^{\infty} \widetilde{c}_{ij} 
\underline{h}_i\right) +  \sum_{j=1}^r \left(\sum_{p=1}^{2n} \widetilde{c}_j 
\epsilon^p \wedge \widetilde{\alpha}_j \otimes  
\sum_{i=1}^{\infty}e_p(\widetilde{c}_{ij}) \underline{h}_i\right) \\
&=& \sum_{i=1}^{\infty} \sum_{j=1}^{r} \left(\widetilde{c}_{ij} 
d_{k-1}(\widetilde{c}_j\widetilde{\alpha}_j) +  \sum_{p=1}^{2n} 
e_p(\widetilde{c}_{ij})\epsilon^p \wedge \widetilde{c}_j \widetilde{\alpha}_j 
\right) \otimes \underline{h}_i\\
&=& \sum_{i=1}^{\infty} \sum_{j=1}^{r} \left(d_{k-1}(\widetilde{c}_{ij}   \widetilde{c}_j   \widetilde{\alpha}_j)\right) \otimes \underline{h}_i
\end{eqnarray*}
since $d_0\widetilde{c}_{ij}=\sum_{p=1}^{2n}e_p(\widetilde{c}_{ij})\epsilon^p.$
Summing-up,  $\sum_{l=1}^sc_{li} \alpha_l = \sum_{j=1}^r 
d_{k-1}(\widetilde{c}_{ij}\widetilde{c}_j\alpha_j)$
for each $i \in \mathbb{N}$.
In particular,  $\sum_{l=1}^s c_{li}\alpha_l$ is exact for each $i\in 
\mathbb{N}.$
Thus, $\sum_{i=1}^{\infty} \sum_{l=1}^s[\alpha_l]\otimes c_{li}h_i = \sum_{i=1}^{\infty} [\sum_{l=1}^sc_{li}\alpha_l]\otimes h_i = 0$ and 
 $[\gamma]$ is injective. 
The verification of the $CH$-equivariance of $[\gamma]$ 
is straightforward.

4) The inverse to $[\gamma]$ which we denote by $[\delta],$ is defined by 
 $[\delta]([\sum_{j=1}^r 
\alpha_j \otimes \sum_{i=1}^{\infty}c_{ij}\underline{h}_i]) = 
\sum_{i=1}^{\infty}[\sum_{j=1}^r
c_{ij}\alpha_j] \otimes h_i$ where $c_{ij}$ are smooth 
functions on $M$ and $\alpha_j \in \Omega^k(M),$ $j=1,\ldots, r$. 
The sum at the right-hand side can be shown to converge by undertaking 
a similar procedure as in the proof of the surjectivity of $[\gamma].$  
Further, it is immediate to check that  $[\delta]$  maps cocycles into 
cocycles and vanishes on coboundaries. Summing-up, $[\delta]$ is well 
defined.
Let us define a map $\Xi(\sum_{j=1}^r\alpha_j \otimes \sum_{i=1}^{\infty}c_{ij} 
\underline{h}_i) = \sum_{i=1}^{\infty}\sum_{j=1}^r c_{ij} \alpha_j \otimes
h_i.$  Since $M$ is compact, the norm of $\Xi$ is easily seen to be finite and 
thus, $\Xi$ is continuous as a linear map between normed spaces.
The continuity of $[\delta]$ is established
by the use of   the   commutative  diagram ($q_1 \otimes \mbox{Id}_H$ and $\Xi$ 
are continuous, and $q_2$ is open)
 $$\begin{xy}\xymatrix{
Z^k_{dR}(M)\otimes H \ar[d]^{ q_1 \otimes \mbox{Id}_H } &  
Z^k(d_{\bullet}^{\Phi^H})  \ar@{->}[l]_-{\Xi} \ar[d]^{q_2} \\ 
H^k_{dR}(M)\otimes H                           &  H^k(d_{\bullet}^{\Phi^H})  
\ar@{->}[l]_-{[\delta]}}
\end{xy}$$
Summing-up, we have
$H^k(d^{\Phi^H}_{\bullet}) \simeq H^k_{dR}(M)\otimes H$ as $CH$-modules (by  item 3) and 
as topological spaces by this item.

5) Notice that we already proved that $H^k(d^{\Phi}_{\bullet})$ and 
$H^k(d^{\Phi^H}_{\bullet})$ are isomorphic   as $CH$-modules and
topological spaces (item 2).
 Since $H^k_{dR}(M)$ is finite dimensional, we may use the projection 
$\mbox{pr}$ defined in item 3. The  quotient topology on the cohomology 
group is generated by the scalar product 
$([\alpha],[\beta])_{H^k_{dR}}=(\mbox{pr}(\alpha),\mbox{pr}(\beta))_{\Omega^k}$ 
where at the right-hand side, the  Hodge scalar product on differential forms 
is meant.
Consequently, $H^k_{dR}(M)\otimes H$  possesses a pre-Hilbert 
$C^*$-product given by $([\alpha]\otimes h, [\beta]\otimes 
h')_{\otimes}=([\alpha],[\beta])_{H^k_{dR}} (h,h')$ where at the 
right-hand side, $(h,h')$ denotes the Hilbert $CH$-product on 
$H$. The norm induced by this
pre-Hilbert 
$CH$-product   coincides with the norm 
induced by the scalar product considered on the tensor product (Lemma 1).
By the finite dimensionality of the de Rham cohomology groups,  
$H^k_{dR}(M)\otimes H$ is a complete inner product space.  
 In particular, the pre-Hilbert $CH$-product on the tensor product of 
the de Rham  cohomology group  and $H$ is a Hilbert $CH$-product.

 Since  $H^k(d_{\bullet}^{\Phi^H})$ is homeomorphic 
and $CH$-isomorphic to $H^k_{dR}(M)  \otimes H$ (item 4)    we can equip it with
the inner product 
$([s_1],[s_2])=\left([\delta]([s_1]), [\delta]([s_2])\right)_{\otimes}$ where 
$s_1,s_2$ are cocycles for $d^{\Phi^H}_{\bullet}$.
Because $[\delta]$ is an isomorphism of $CH$-modules, the resulting 
map is  a well defined 
Hilbert $C^*$-product. Since $[\delta]$  is a homeomorphism,
 the topology on 
$H^k(d_{\bullet}^{\Phi^H})$ generated by the norm induced by the scalar product 
 coincides with the 
quotient topology on $H^k(d_{\bullet}^{\Phi^H}).$ Similarly, by 
the use of the map $j$ (item 1 and item 2),
 we induce the Hilbert $CH$-product from $H^k(d_{\bullet}^{\Phi^H})$ to
 $H^k(d_{\bullet}^{\Phi})$ making it a Hilbert $CH$-module isomorphic to 
$H^k_{dR}(M) \otimes H$ as well.
 
Because  $H^k_{dR}(M)$ is finite dimensional and 
 $H$ is finitely generated (Lemma 2),  $H^k_{dR}(M)\otimes H$ 
is also finitely generated.
Consequently, it is projective by the Magajna theorem.
Summing-up, $H^k(d^{\Phi}_{\bullet})$  are finitely generated projective Hilbert
$CH$-modules for each $k.$
 \hfill\(\Box\)

\bigskip

We give a much shorter proof of Thm. 18 based on the Hodge theory
for $CH$-bundles.

\bigskip

{\it Second proof of Theorem 18.}
We use the Mishchenko--Fomenko theory elaborated for complexes in  
\cite{KryslAGAG2} and \cite{KryslJGP2}.
It is immediate to compute that the symbols $s_k$ of $d_k^{\Phi}$ are given 
by
$s_k(\alpha \otimes s, \tau) = (\tau \wedge \alpha) \otimes s,$ 
$\tau \in T_m^*M,$ $\alpha \otimes s \in (\mathcal{H}_R^k)_m$ (see 
\cite{KryslDGA}),
and a matter of multilinear algebra to see that they form an  exact sequence in 
places $k=0,\ldots, 2n$ for any 
$\tau \neq 0.$ In the places $k=0,\ldots, 2n,$ the Kuiper complex is thus 
elliptic. Due to Lemma 2, it is a complex on finitely generated projective
Hilbert $CH$-bundles.
Consequently (compactness and ellipticity), its cohomology groups are finitely 
generated and projective 
Hilbert $CH$-modules by Thm. 9 in Kr\'ysl  \cite{KryslJGP2}.   In places 
$\mathbb{Z}\setminus \{0,\ldots, 2n\},$ the cohomology groups
are zero. 
\hfill\(\Box\)

\bigskip

{\bf Remark:} 1) Since the cohomology groups are Hilbert $C^*$-modules, 
the images of $d_{\bullet}^{\Phi}$ and $d_{\bullet}^{\Phi^H}$ are necessarily 
closed.

2) The Hodge theory holds for the Kuiper complex by  Thm. 9 in \cite{KryslJGP2} as 
well. From the proof of Theorem 18 (and the proof of Lemma 2), we see  that 
the $CH$-rank of $H^k(d_{\bullet}^{\Phi})$ equals to the $k$-th Betti number
of $M.$ (See Baki\'c, Gulja\v s \cite{Bakic} for a definition of the $CH$-rank 
if needed.)

\bibliographystyle{amsplain}

\begin{thebibliography}{99}


\bibitem{ALMP} Albin, P., Leichtnam, E., Mazzeo, R., Piazza, P., Hodge theory 
on Cheeger spaces,   \url{https://doi.org/10.1515/crelle-2015-0095}, on-line 
published in 2016. ArXiv--preprint \url{https://arxiv.org/abs/1307.5473}. 

\bibitem{Bakic} Baki\'c, D., Gulja\v s, B., Operators on Hilbert $H^*$-modules, J. Oper. Theory 46 (2001), No. 1,
123--137.

\bibitem{BirkhoffNeumann} Birkhoff, G., von Neumann, J.,  The logic of 
quantum mechanics,  Ann. Math. (2) 37, Princeton (1936), 823--843.

\bibitem{Bong} Bongioanni, B., Torrea, J., 
Sobolev spaces associated to the harmonic oscillator. In: Proc. Indian Acad. Sci. Math. Sci. 116 (2006), No. 3, 337--360.

\bibitem{BorelWallach} Borel, A., Wallach, N., Continuous cohomology, discrete 
subgroups, and representations of reductive groups. Second edition. Mathematical 
Surveys and Monographs, 67. American Mathematical Society, Providence, RI, 2000.

\bibitem{Bruhat} Bruhat, F., Sur les repr\'esentations induites des groupes de Lie, 
Bull. Soc. Math. France 84, 1956, 97--205.

\bibitem{CGLR} Cahen, M., Gutt, S., La Fuente Gravy, L., Rawnsley, J., On 
$Mp^c$-structures and symplectic Dirac operators. J. Geom. Phys. 86 (2014), 
434--466. 

\bibitem{Connes} Connes, A., Noncommutative geometry. Academic Press, Inc., San 
Diego, CA, 1994. 


 \bibitem{DixmierDouady} Dixmier, J., Douady, A.,
Champs continus d'espaces hilbertiens et de $C^*$-alg\'ebres, 
Bull. Soc. Math. France 91 (1963), 227--284. 

\bibitem{Dixmier} Dixmier, J.,  Les $C^*$-alg\'ebres et leurs 
repr\'esentations. Gauthier-Villars, Paris, 1969. 


\bibitem{DM} Dubois-Violette, M., Michor, P., 
Connections on central bimodules in noncommutative differential geometry, J. 
Geom. Phys. 20 (1996), Issue 2--3, 218--232.



\bibitem{FathiGabr} Fathizadeh, F., Gabriel, O., On the 
Chern--Gauss--Bonnet theorem and conformally twisted spectral triples for 
$C^*$-dynamical systems, SIGMA Symmetry Integrability Geom. Methods Appl. 12 
(2016), paper No. 016.


\bibitem{Folland} Folland, G., Harmonic analysis in phase space. Annals of 
Mathematics Studies, 122. Princeton Univ. Press, Princeton, NJ, 1989. 


\bibitem{FM} Fomenko, A., Mishchenko, A., The index of elliptic operators over 
$C^*$-algebras, Math. USSR-Izv. 15, (1980), 87--112.


\bibitem{FH} Forger, M., Hess, H., Universal metaplectic structures and 
geometric quantization, Comm. Math. Phys. 64 (1979), No. 3, 269--278. 


\bibitem{Freed}  Freed, D., Lott, J., An index theorem in differential 
$K$-theory, Geom. Topol. 14 (2010), No. 2, 
903--966. 

\bibitem{Goodman} Goodman, R., Analytic and entire vectors for representations 
of Lie groups. Trans. Amer. Math. Soc. 143 (1969), 55--76.

\bibitem{Groth} Grothendieck, A., Produits tensoriels topologiques et Espaces 
nucl\'eaires. Amer. Mat. Soc., Providence, RI, 1966.

\bibitem{HabermannAn}  Habermann, K., The Dirac operator on symplectic spinors.
Annals Glob. Anal. Geom. 13 (1995), No. 2, 155--168. 

\bibitem{HH} Habermann, K., Habermann, L.,  Introduction to Symplectic Dirac Operators, Lecture Notes in Mathematics 1887.
Springer-Verlag, Berlin, 2006.  

\bibitem{Hain} Hain, R., The Hodge--de Rham theory of modular groups, In: Recent 
advances in Hodge theory, London Math. Soc. Lecture Note Ser., 427. 
Cambridge Univ. Press, Cambridge, 2016, 422--514.

\bibitem{HenTeitel} Henneaux, M.,  Teitelboim, C., Quantization of Gauge 
Systems, Princeton University Press, Princeton, NJ, 1992. 

\bibitem{HerczegWaldron} Herczek, G.,  Waldron, A., Contact geometry and 
Quantum mechanics,  Physics Letters, Section B: Nuclear, Elementary Particle and 
High-Energy Physics, 781, 2018, 312--315, 
\url{https://doi.org/10.1016/j.physletb.2018.04.008}.

\bibitem{Hodge} Hodge, W.,  The theory and applications of harmonic integrals. 
2nd ed. Cambridge,  at the Univ. Press, 1952.

\bibitem{Illusie}  Illusie, L., Complexes quasi-acycliques directs de fibr\'es 
banachiques, C. R. Acad. Sci. Paris 260, 1965, 6499--6502.

\bibitem{Keyl} Keyl, M., Kiukas, J., Werner, R.,
Schwartz operators, Rev. Math. Phys. 28 (2016), No. 3, 1630001, 60 pp.

\bibitem{Kim} Kim, J., A splitting theorem for holomorphic Banach bundles, 
Math. Z. 263, 2009, No. 1, 89--102.

\bibitem{Kirillov} Kirillov, A., Lectures on the orbit method. Graduate 
Studies in Mathematics, 64. American Mathematical Society, Providence, RI, 2004. 
ISBN 0-8218-3530-0. 

\bibitem{Klein} Klein, F., \"Uber die Entwicklung der Mathematik im 19. 
Jahrhundert. Band 1, Springer, Berlin, 1926. Engl. transl. by M. Ackermann: 
Development of Mathematics in the 19th Century with Appendices 
"Kleinian Mathematics from an Advanced Standpoint" by  R. Hermann, 1st Edition, 
Math Sci Press, Brookline, MA, 1979.

\bibitem{KnappOverview} Knapp, A., Representation Theory of Semisimple Groups. An Overview Based on Examples, 
Reprint of the 1986 original. Princeton Landmarks in Mathematics, Princeton Univ. Press, Princeton, 2001.

\bibitem{KMS} Kol\'a\v{r}, I., Michor, P., Slov\'ak, J., Natural operations in 
differential geometry. Springer-Verlag, Berlin, 1993. 

\bibitem{Kostant} Kostant, B., Symplectic spinors, In: Symposia Mathematica, Vol. XIV 
(Convegno di Geometria Simplettica e Fisica Matematica, INDAM, Rome, 1973). Academic Press, London, 1974, 139--152.

\bibitem{KryslDGA0} Kr\'ysl, S., Classification of 1st order symplectic spinor operators over contact projective geometries,
Vol. 26 (2008) Issue 5,  553--565.

\bibitem{KryslDGA} Kr\'ysl, S., Cohomology of the de Rham complex twisted by the 
oscillatory representation, Diff. Geom. Appl. 33 (2014), No. 5,   
 290--297.

 \bibitem{KryslAGAG2} Kr\'ysl, S., Hodge theory for self-adjoint parametrix 
possessing complexes over $C^*$-algebras, Annals Glob. Anal. Geom. 47 
 (2015), No. 4, 359--372.

\bibitem{KryslJGP2} Kr\'ysl, S., Elliptic complexes over $C^*$-algebras of 
compact operators, Journal of Geometry and Physics 101 (2016), 27--37.

\bibitem{KryslHabil} Kr\'ysl, S., Hodge Theory and Symplectic Spinors. 
Habilitation thesis. Faculty of Mathematics and Physics, Charles University, 
Prague, 2016, electronically available at
\url{http://msekce.mff.cuni.cz/~krysl/habil.pdf},
short version at ArXiv \url{https://arxiv.org/abs/1708.02026}.
 
\bibitem{Kuiper} Kuiper, N., The homotopy type of the unitary group of Hilbert space, Topology 3 (1965), 19--30. 

\bibitem{Larrain} Larra\'in-Hubach, A., $K$-theories for classes of infinite 
rank bundles. In: Analysis, geometry and quantum field theory,
Contemp. Math., 584. Amer. Math. Soc., Providence, RI, 2012, 79--97. 

\bibitem{Lempert} Lempert, L., On the cohomology groups of holomorphic Banach 
bundles, Trans. Amer. Math. Soc. 361 (2009), No. 8, 4013--4025. 


\bibitem{Li} Li, P., Cauchy--Schwarz-type inequalities on K\"{a}hler 
manifolds--II, J. Geom. Phys. 99 (2016), 256--262. 

\bibitem{Ludwig} Ludwig, G., Versuch einer axiomatischen Grundlegung der 
Quantenmechanik und allgemeinerer physikalischer Theorien, Zeitschrifft f\"ur 
Physik 181 (1964), Issue 3,  233--260.

\bibitem{MaedaRosenberg} Maeda, Y., Rosenberg, S., Traces and characteristic 
classes in infinite dimensions, Geometry and analysis on manifolds, 413--435,
Progr. Math., 308, Birkh\"{a}user/Springer, Cham (CH), 2015. 

\bibitem{Magajna} Magajna, B., Hilbert $C^*$-modules in which all closed submodules are complemented,
 Proc. Amer. Math. Soc. 125 (1997), No.  3,  849--852.

\bibitem{MMS} Mathai, V., Melrose, R., Singer, I., The index of projective 
families of elliptic operators: the decomposable case. In:
Dai, X. (ed.) et al., From probability to geometry II (Volume in honor of 
the 60th birthday of Jean-Michel Bismut). Paris SMF,
Ast\'erisque 328, 255--296 (2009).

\bibitem{Maurin} Maurin, K.,  The Riemann legacy. Riemannian ideas in 
mathematics and physics. Mathematics and its Applications, 417. Kluwer Academic 
Publishers Group, Dordrecht, 1997.

\bibitem{Michor} Michor, P.,  Topics in differential geometry.  
Graduate Studies in Mathematics 93. Providence, RI, American Mathematical Society, 2008. 

\bibitem{Mishchenko}  Mishchenko, A., The theory of elliptic operators over 
$C^*$-algebras, Dokl. Akad. Nauk SSSR 239 (1978), No. 6, 1289--1291. 

\bibitem{Neeb} Neeb, K.,  On differentiable vectors for representations of 
infinite dimensional Lie groups. J. Funct. Anal. 259 (2010), No. 11, 2814--2855.

\bibitem{Nekovar}  Nekov\'a\v{r}, J., Scholl, A., Introduction to plectic 
cohomology. In: Jiang, D., Shahidi, F., Soudry, D. (eds.), Advances in the theory of automorphic forms and their 
$L$-functions, 321--337, Contemp. Math., 664, Amer. Math. Soc., Providence, RI, 2016. 

\bibitem{NeumannKM} von Neumann, J., Zur Operathorenmethode in der klassischen 
Mechanik, Ann. Math. (2) 33 (1932), 257--642. 

\bibitem{Palais}  Palais, R., Seminar on the Atiyah--Singer index 
theorem (with contributions by M. F. Atiyah, A. Borel, E. E. Floyd, R. T. 
Seeley, W. Shih and R. Solovay), Annals of Mathematics Studies, No. 57. Princeton 
Univ. Press, Princeton, NJ, 1965.

\bibitem{Paschke} Paschke, W., Inner product modules over $B^*$-algebras, 
Trans. Amer. Math. Soc. 182 (1973),  443--468.

\bibitem{Raeburn} Raeburn, I., Williams, D., Morita equivalence and continuous-trace $C^*$-algebras. 
Mathematical Surveys and Monographs, 60. American Mathematical Society, Providence, RI, 1998. 
ISBN 0-8218-0860-5. 

\bibitem{Rieffel}   Rieffel, M., Induced representations of
$C^*$-algebras, Advances in Mathematics 13 (2) (1973),
176--257.

\bibitem{RR} Robinson, P.,  Rawnsley, J.,  The metaplectic representation, $Mp^c$-structures and geometric quantization, 
Mem. Am. Math. Soc. 410, 1989.

\bibitem{Rohrl} R\"{o}hrl, H., \"{U}ber die Kohomologie berechenbarer Fr\'echet 
Garben, Comment. Math. Univ. Carolinae 10 (1969), 625--640. 

\bibitem{Rordam} R\o rdam, M., Larsen, F., Laustsen, N. J., An introduction to $K$-theory for $C^*$-algebras, London Mathematical Society, Cambridge University Press,
 2000.

\bibitem{Rosenberg} Rosenberg, J., Topology, $C^*$-algebras, and string duality.
CBMS Regional Conference Series in Mathematics, 111. Published for the Conference Board of the Mathematical Sciences, Washington, DC,
American Mathematical Society, Providence, RI, 2009.

\bibitem{Ryan} Ryan, R., Introduction to tensor products of Banach Spaces,  
Springer Monographs in Mathematics.  London, 2002.

\bibitem{Schick} Schick, T., $L^2$-index theorems, $KK$-theory, 
and connections, New York J. Math. 11 (2005), 387--443.  
 
\bibitem{SchmidVilonen} Schmid, W., Vilonen, K., Hodge theory and 
unitary representations. Representations of reductive groups, 443--453, Progr. 
Math., 312, Birkh\"{a}user/Springer, Cham (CH), 2015. 

\bibitem{Schotten} Schottenloher, M., The unitary group and its strong 
topology. ArXiv--preprint, \url{https://arxiv.org/abs/1309.5891}.

\bibitem{Shale} Shale, D., Linear symmetries of free boson fields, Trans. Amer. Math. Soc. 103 (1962),  149--167.

\bibitem{ST} Solovyev (Solov'\"{e}v), Y., Troitsky (Troickij), E., $C^*$-algebras and 
elliptic operators in differential topology. Transl. of Mathem. Monographs 
192,  Amer. Math. Soc., Providence, RI, 2001.

\bibitem{Troitsky} Troitsky, E., The index of equivariant elliptic operators 
over $C^*$-algebras, Annals Glob. Anal. Geom.,  5 (1987), No.  1,  3--22.

\bibitem{Tsai} Tsai, C., Tseng, L., Yau, S., Cohomology and 
Hodge theory on symplectic manifolds: III. J. Differential Geom. 103 (2016), No. 
1, 83--143.

\bibitem{Wallach} Wallach, N.,  Symplectic geometry and Fourier analysis. With 
an appendix on quantum mechanics by Robert Hermann. Lie Groups: History, 
Frontiers and Applications, Vol. V. Math Sci Press, Brookline, MA, 1977.

\bibitem{Warner} Warner, G.,  Harmonic analysis on semi-simple Lie groups I. 
Die Grundlehren der mathematischen Wissenschaften, Band 188. Springer-Verlag, 
New York--Heidelberg, 1972. 

\bibitem{Weil} Weil, A., Sur certains groupes d'op\'erateurs unitaires, Acta Math. No.  111 (1964),  143--211. 

\end{thebibliography}

\end{document}